\documentclass[preprint,12pt]{elsarticle}
\UseRawInputEncoding
\usepackage{amssymb}
\usepackage{amsmath}
\usepackage{geometry}
\allowdisplaybreaks[4]

\journal{......}

\begin{document}

\begin{frontmatter}

%% Title, authors and addresses

%% use the tnoteref command within \title for footnotes;
%% use the tnotetext command for the associated footnote;
%% use the fnref command within \author or \address for footnotes;
%% use the fntext command for the associated footnote;
%% use the corref command within \author for corresponding author footnotes;
%% use the cortext command for the associated footnote;
%% use the ead command for the email address,
%% and the form \ead[url] for the home page:
%%
%% \title{Title\tnoteref{label1}}
%% \tnotetext[label1]{}
%% \author{Name\corref{cor1}\fnref{label2}}
%% \ead{email address}
%% \ead[url]{home page}
%% \fntext[label2]{}
%% \cortext[cor1]{}
%% \address{Address\fnref{label3}}
%% \fntext[label3]{}

\title{Global existence of large solutions to 3-D incompressible Navier-Stokes system}

%% use optional labels to link authors explicitly to addresses:
%% \author[label1,label2]{<author name>}
%% \address[label1]{<address>}
%% \address[label2]{<address>}

\author{Shaolei Ru\\}

\address{~~~~~~School of Mathematical Sciences, Zhejiang Normal University, China~~~~~~~~\\
}

\begin{abstract}
%% Text of abstract
This paper proves that the 3-D Navier-Stokes system has a unique global solution under an assumpution on
the initial data. That allow the data to be arbitrarily large in the scale invariant space $\dot{B}_{\infty,\infty}^{-1}$, which contains all the known
spaces in which there is a global solution for small data.
Moreover, this result provides various examples of arbitrarily large initial data giving rise to global solutions.
\end{abstract}

\begin{keyword}
%% keywords here, in the form: keyword \sep keyword
Navier-Stokes equations; Global well-posedness; Large initial data.
%% MSC codes here, in the form: \MSC code \sep code
%% or \MSC[2008] code \sep code (2000 is the default)
\MSC[2020] 35Q30, 76D03.\\
E-mail address:~rushaolei@zjnu.edu.cn
\end{keyword}

\end{frontmatter}

%%
%% Start line numbering here if you want
%%
% \linenumbers

%% main text
\section{Introduction}

 This paper considers the 3D incompressible Navier-Stokes (NS) equations on $\mathbb{R}^{3}$ read:
\begin{align}\label{hang113}(NS)\left\{
  \begin{array}{lll}
  \partial_{t}u-\Delta u+(u\cdot \nabla)u+\nabla p=0,\\
  div u=0,~u(0,x)=u^{0}(x),\\
  \end{array}
  \right.
\end{align}
where $\Delta=\sum_{i=1}^{3}\partial_{i}^{2}$,
$\nabla=(\partial_{1},\partial_{2},\partial_{3})$,
$div
u=\partial_{1}u_{1}+\partial_{2}u_{2}+\partial_{3}u_{3}$
and $u=(u_{1},u_{2},u_{3})$. Here $u:[0,T)\times
\mathbb{R}^{3}\rightarrow\mathbb{R}^{3}$ stands for the velocity
field and $p:[0,T)\times \mathbb{R}^{3}\rightarrow\mathbb{R}^{3}$
for the pressure. The kinematic viscosity of the
fluid has been chosen equal to one for simplicity. Note that Eq.(\ref{hang113}) is scaling invariant in the following sense: if the pair $(u,p)$ is a solution of (\ref{hang113}) on the time
interval [0, T), then $(u_{\lambda},p_{\lambda})$ defined by
\begin{align}\label{18}u_{\lambda}(t,x)=\lambda u(\lambda^{2}t,\lambda x)~~ and~~
p_{\lambda}(t,x)=\lambda^{2} u(\lambda^{2}t,\lambda x)\end{align} is a solution on the time interval $[0,\lambda^{2}T)$ with initial data
$u_{\lambda}^{0}=\lambda u^{0}(\lambda x)$. A function space $X$ defined in $\mathbb{R}^{d}$ is said to be a critical space for Eq.(\ref{hang113}) if the norms of
$u_{\lambda}^{0}$ in $X$ are equivalent for all $\lambda>0$.

Let $(u,p)$ be a smooth solution of the NS equation. Taking the
divergence of the first equation in (NS) and noticing that $div
u=0$, one can immediately obtain that
$$\Delta p+div[(u\cdot\nabla)u]=0.$$
It follows that $\nabla p=(-\Delta)^{-1}\nabla
div[(u\cdot\nabla)u]$. For convenience, denote
$$\mathbb{P}:=I+(-\Delta)^{-1}\nabla div.$$ Inserting $\nabla p$ into Eq.(NS), we get
$$\partial_{t}u+\mathbb{P} [(u\cdot\nabla) u]-\Delta u=0,~u(0,x)=u_{0}(x).$$
Based on the above observation, one sees that the NS equation belongs to a nonlinear parabolic
equation. Besides, noticing that $divu=0$, we have $(u\cdot \nabla)u=div(u\otimes u)$.

A huge literature study the well-posedness of the NS and related equations, for example, \cite{giga14}\cite{Coron}, etc. Leray \cite{Leray} initiated
the study of the NS equations by proving the existence
of weak solutions. The study of the equations for small data
started with a groundbreaking paper by Fujita and Kato \cite{Fujita}, in which they
showed local well-posedness and small data global well-posedness in the space $\dot{H}^{\frac{1}{2}}(\mathbb{R}^{3})$.
Later, Kato \cite{kato12}
 studied the NS equations by proving that the NS
problem (\ref{hang113}) is locally well-posed in $L^{3}(\mathbb{R}^{3})$ and global well-posed
if the initial data are small in $L^{3}(\mathbb{R}^{3})$. The method consists in applying a Banach fixed point theorem to the
integral formulation of the equation, and was generalized by M. Cannone, Y. Meyer and F.
Planchon in \cite{cannone} to Besov spaces of negative index of regularity. More precisely they proved
that if the initial data is small in the Besov space $\dot{B}_{p,\infty}^{-1+\frac{3}{p}}(\mathbb{R}^{3})$ (for $p<\infty$), then there is a unique,
global in time solution. It is pointed out that this result allows to construct global solutions for strongly oscillating initial data which may have a large norm in $\dot{H}^{\frac{1}{2}}(\mathbb{R}^{3})$ or in $L^{3}(\mathbb{R}^{3})$. A typical example is \begin{eqnarray}\label{A21}u_{\varepsilon}^{0}(x)=\varepsilon^{-\alpha}sin(\frac{x_{3}}{\varepsilon})(-\partial_{2}\varphi(x),\partial_{1}\varphi(x),0),
\end{eqnarray}
where $0<\alpha<1$ and $\varphi\in \mathcal{S}(\mathbb{R}^{3};\mathbb{R})$. More recently in \cite{koch17}, H. Koch and D. Tataru obtained a unique global in time solution for data small enough in a more general space, consisting of vector fields whose components are derivatives of BMO functions. The norm in that space is given by
$$\|u^{0}\|_{BMO^{-1}}^{2}=sup_{t>0}t\|e^{t\Delta}u^{0}\|_{L^{\infty}}^{2}+sup_{x\in \mathbb{R}^{d}, R>0}\frac{1}{R^{d}}\int_{P(x,R)}|(e^{t\Delta}u_{0})(t,y)|^{2}dy,$$where $P(x,R)$ stands for the parabolic set $[0,R^{2}]\times B(x,R)$.
Whereafter, Bourgain and Pavlovi\'{e}
\cite{bourgain23} proved that the Cauchy problem for the three-dimensional Navier-Stokes equations is ill-posed in
$\dot{B}_{\infty,\infty}^{-1}(\mathbb{R}^{3})$ in the sense that a "norm inflation" happens in finite time.

 It is easy to see that the following spaces are critical spaces for 3D NS:
$$\dot{H}^{\frac{1}{2}}(\mathbb{R}^{3})\hookrightarrow L^{3}(\mathbb{R}^{3})\hookrightarrow \dot{B}_{p,\infty}^{-1+\frac{3}{p}}(\mathbb{R}^{3})\hookrightarrow BMO^{-1}(\mathbb{R}^{3})\hookrightarrow \dot{B}_{\infty,\infty}^{-1}(\mathbb{R}^{3}).$$
 This observation shows that if one wants to go beyond a smallness assumption on the initial data to prove the global existence of unique solutions, one should check that the $\dot{B}_{\infty,\infty}^{-1}(\mathbb{R}^{3})$ norm of the initial data may be chosen large. Unfortunately, the example stated in (\ref{A21}) can not  achieve these ends. Some progress has been made on this question. The authors in \cite{Chemin1} constructed an example of periodic initial data which is large in $B_{\infty,\infty}^{-1}(\mathbb{T}^{3})$ but yet generates a global solution. Such an initial data is given by
\begin{eqnarray}\label{A22}u_{N}^{0}=(N\nu_{h}^{0}(x_{h})cos(Nx_{3}),-div_{h}\nu_{h}^{0}(x_{h})sin(Nx_{3})),\end{eqnarray}
where $x_{h}=(x_{1},x_{2})$ and $\|\nu_{h}^{0}(x_{h})\|_{L^{2}(\mathbb{T}^{2})}\leqslant C(ln N)^{\frac{1}{4}}$.
  This was generalized to the case of the space $\mathbb{R}^{3}$ in \cite{Chemin2}. Thereafter, in \cite{Chemin3} the authors showed a global
existence result for large data which are slowly varying in one direction. More
precisely, they proved that for $\varepsilon$ small enough, the initial data
\begin{eqnarray}\label{A23}u_{\varepsilon}^{0}(x_{h},x_{3})=(\nu_{h}^{0}(x_{h},\varepsilon x_{3}),0)
+(\varepsilon\omega_{h}^{0}(x_{h},\varepsilon x_{3}),\omega_{3}^{0}(x_{h},\varepsilon x_{3})),\end{eqnarray}
where $(\nu_{h}^{0},0)$ and $\omega^{0}$ are two smooth divergence free vector fields, generates
a global smooth solution. In this case, its $\dot{B}_{\infty,\infty}^{-1}(\mathbb{R}^{3})$
norm is large but does not depend on the slow parameter. Subsequently, J. Chemin, I. Gallagher and M. Paicu \cite{Chemin4} proved a global
existence result for the initial data
$$u_{0,\varepsilon}=(\nu_{h}^{0}(x_{h},\varepsilon x_{3}),\frac{1}{\varepsilon}\nu_{3}^{0}(x_{h},\varepsilon x_{3}))$$
on $\mathbb{T}^{2}\times\mathbb{R}$. Where $div \nu^{0}=0$ and $\|e^{aD_{3}}\nu^{0}\|_{H^{4}}$ sufficiently small.
In this situation, the $B_{\infty,\infty}^{-1}(\mathbb{T}^{2}\times\mathbb{R})$ norm blows up
as the small parameter goes to zero. In the case of the whole space $\mathbb{R}^{3}$, it seems
difficult to control the very low horizontal frequencies.
M. Paicu and Z. Zhang considered in \cite{Paicu} an intermediate situation on $\mathbb{R}^{3}$ between the above two cases later on. For some other results, one can refer to \cite{Foias}, etc.

On the other hand, some papers studied the condition one can assume on the
initial data for the Eq.(1) to be globally wellposed. In \cite{Chemin1}, the author gave a nonlinear smallness condition on the initial data for the periodic, three dimensional, incompressible Navier-Stokes
equations to be globally wellposed (for the details, see Theorem D in Appendix). It is not a smallness condition, as the data is allowed to be arbitrarily large in the scale invariant space $\dot{B}_{\infty,\infty}^{-1}$. Recently, in \cite{Liu} the authors proved that the Eq.(\ref{hang113}) has a unique global Fujita-Kato
solution provided that the $H^{-\frac{1}{2},0}$ norm of $\partial_{3}u^{0}$ is sufficiently small compared to $exp\{-C(A_{\delta}(u_{h}^{0})+B_{\delta}(u^{0}))\}$, with $A_{\delta}(u_{h}^{0})$ and $B_{\delta}(u^{0})$ being scaling invariant quantities of the initial data (for the details, see Theorem E in Appendix). This result also provides some classes of large initial data which generate unique global solutions to Eq.(\ref{hang113}). For some other results, one can refer to \cite{Liu1}, etc.

The key ingredient used in the above is that one can decompose the solution of (NS) as a sum of a large two-dimensional
solution of (NS) and a small three-dimensional one. This article is not in this situation.

Inspired by the statement in \cite{Fefferman} and \cite{Chemin4}, this paper mainly study the condition and size of the initial data yielding global existence of solutions
to the Eq.(\ref{hang113}). The first purpose of this text is to establish a diverse condition for regular initial data for the Eq.(\ref{hang113}) to be globally wellposed. In the proof, we need to use a different phase function, which based on a decomposition of the solution and the special structure of the equations.
The main idea of the approach is that if $u$ denotes the solution of (NS) associated with $u^{0}$, which exists at least for a short time since $u^{0}\in \dot{B}_{2,\gamma}^{\frac{1}{2\gamma^{\prime}},\frac{1}{2\gamma}}$, then it can be decomposed as follows:$$u(t)=\bar{u}(t)+\nu(t),$$ where $\bar{u}(t)$ satisfies the following equation:
\begin{align}\label{hang111}\left\{
  \begin{array}{lll}
  \partial_{t}\bar{u}_{3}-\Delta \bar{u}_{3}+(u\cdot \nabla)u_{3}=0,\\
  \partial_{t}\bar{u}_{i}-\Delta \bar{u}_{i}+div_{h}(\bar{u}_{h}\bar{u}_{i})=0,~~i=1,2,\\
  \bar{u}(0,x)=u^{0}(x).\\
  \end{array}
  \right.
\end{align}
One notices that the vector field $\nu(t)$ satisfies the perturbed Navier-Stokes system:
\begin{align}\label{hang112}\left\{
  \begin{array}{lll}
  \partial_{t}\nu_{3}-\Delta \nu_{3}+[\mathbb{Q}(u\cdot \nabla)u]_{3}=0,\\
  \partial_{t}\nu_{i}-\Delta \nu_{i}+div_{h}(\nu_{h}\bar{u}_{i})+div_{h}(u_{h}\nu_{i})+\partial_{3}(u_{3}u_{i})+[\mathbb{Q}(u\cdot \nabla)u]_{i}=0,\\
  \nu(0,x)=0,~~i=1,2,\\
  \end{array}
  \right.
\end{align}
where $[\mathbb{Q}(u\cdot \nabla)u]_{j}=R_{j}\sum_{i=1}^{3}R_{i}[(u\cdot \nabla)u_{i}]$ and $R_{j}$, j=1,2,3, is the Riesz transform defined by $$\widehat{R_{j}f}(\xi)=-i\frac{\xi_{j}}{|\xi|}\widehat{f}(\xi).$$ The proof of main results consists in studying both systems. Through observation, from the construction of equation (\ref{hang111}), one may study some special estimates for horizontal components of the solution. Using a different phase function, one can further study a priori estimates for $\bar{u}$. On the other hand, in view of the special structure of $\mathbb{Q}[(u\cdot \nabla)u]$ and its relationship to the divergence of $\bar{u}$ in Eq.(\ref{hang111}), one can get the global well-posedness of the perturbed system (\ref{hang112}) with small data.
As the second purpose of this text, this technology allow us to consider more wider classes of arbitrarily large initial data giving rise to global solutions.
In particular, that allow the components of initial data, both $u_{2}^{0}$ and $u_{3}^{0}$, to be arbitrarily large in $\dot{B}_{\infty,\infty}^{-1}(\mathbb{R}^{3})$.\\\\
 \begin{bf}Main results:\end{bf}\\\\
\begin{bf}Theorem 1.\end{bf} Let $a\geqslant1$, $1<\gamma<2$ and $divu^{0}=0$. Then there exist constants $C, \lambda>0$ and a small enough positive constant $\eta$ such that if
\begin{align}\label{hang3332}&C\exp\{8C\lambda\|e^{aD_{3}}u_{2}^{0}\|_{\dot{B}_{2,\gamma}^{\frac{1}{2\gamma^{\prime}},\frac{1}{2\gamma}}}^{\gamma^{\prime}}\}(\|e^{aD_{3}}u^{0}\|_{\dot{B}_{2,\gamma}^{\frac{1}{2\gamma^{\prime}}-\frac{1}{\gamma},\frac{1}{2\gamma}}}+\|e^{aD_{3}}u_{1}^{0}\|_{\dot{B}_{2,\gamma}^{\frac{1}{2\gamma^{\prime}},\frac{1}{2\gamma}}}\\
&~~~~~~~~~~~~~~~~~~~~~~~~~~~~~~~~~~~~~~~~~~~~~~~~~~~~~~
+\|e^{aD_{3}}\partial_{3}u_{3}^{0}\|_{\dot{B}_{2,\gamma}^{\frac{1}{2\gamma^{\prime}}-1,\frac{1}{2\gamma}}})<\eta,\nonumber\end{align} the Navier-Stokes system has a unique global solution on $\mathbb{R}^{3}$.

Using the results in Theorem 1 and scaling invariance, one can establish the following global well-posedness result for (NS) with more wider classes
of large initial data:\\\\
\begin{bf}Theorem 2.\end{bf} Let $0\leqslant\alpha,\beta,\sigma<1$.  Assume that $e^{aD_{3}}\omega^{0}\in B_{1,1}^{2}$, be a divergence free vector field on $\mathbb{R}^{3}$. There exists a positive constant $\varepsilon_{0}$ such that if $\varepsilon<\varepsilon_{0}$, then the initial data
\begin{align*}&u_{\varepsilon,\alpha,\beta,\sigma}^{0}(x)\\
&=(\varepsilon^{1-\alpha-\sigma} \omega_{1}^{0}(\varepsilon^{\alpha}x_{1},\varepsilon^{\beta} x_{2},\varepsilon x_{3}),\varepsilon^{1-\beta-\sigma} \omega_{2}^{0}(\varepsilon^{\alpha}x_{1},\varepsilon^{\beta} x_{2},\varepsilon x_{3}),\varepsilon^{-\sigma}\omega_{3}^{0}(\varepsilon^{\alpha}x_{1},\varepsilon^{\beta} x_{2},\varepsilon x_{3})),\end{align*} with $1-\frac{3}{2}\alpha-\frac{1}{2}\beta>\sigma$ and $1-\frac{1}{2}\alpha-\frac{3}{2}\beta>\sigma$,
generates a unique global smooth solution of (NS).\\\\
\begin{bf}Remarks.\end{bf}
\begin{itemize}\setlength{\itemsep}{0pt}\setlength{\parsep}{0pt}\setlength{\parskip}{0pt}
\item The definition of the anisotropic functional spaces appeared in the above Theorems will be given in Section 2. Comparing to the previous results in \cite{Chemin1} and \cite{Liu}, etc, the condition appearing in Theorem 1 is much simpler and allow us to consider the global well-posedness result for (NS) with more wider classes of large initial data.
 \item From the statement in \cite{Chemin3} and subsection 6.3, one can easily see that $$\|u_{\varepsilon,\alpha,\beta,\sigma}^{0}(x)\|_{\dot{B}_{\infty,\infty}^{-1}}\thicksim \varepsilon^{-(\alpha\wedge\beta)-\sigma}.$$ While $\alpha,\beta,\sigma=0$, this is the case in \cite{Paicu3} and \cite{Liu}. When $\alpha,\beta=0$ and $\sigma=1$, this is the case studied
by \cite{Chemin4}.  For the later case, the data was supposed to evolve
in a special domain $\mathbb{T}^{2}\times\mathbb{R}$. The intermediate case $\alpha,\beta=0, \sigma=\frac{1}{2}$ and $\alpha,\beta=0, \frac{1}{2}<\sigma<1$ was studied respectively by \cite{Paicu5} and \cite{Paicu}. The data, for $0<\alpha,\beta<1$, is not contained in the previous results. In this case, both $u_{2}^{0}$ and $u_{3}^{0}$ can be arbitrarily large in $\dot{B}_{\infty,\infty}^{-1}$.
\item Similar to the arguments as in subsection 6.3, we can also consider a
different type of initial data, with larger amplitude but strongly oscillating in the
horizontal variables, namely, initial data of the form
\begin{align*}\widetilde{u}_{\varepsilon}^{0}(x)&=(0,\varepsilon^{-(\beta+\sigma)} \omega_{2}^{0}(\varepsilon^{\alpha-1}x_{1},\varepsilon^{\beta-1} x_{2},x_{3}),\varepsilon^{-(1+\sigma)}\omega_{3}^{0}(\varepsilon^{\alpha-1}x_{1},\varepsilon^{\beta-1} x_{2},x_{3})),\end{align*}
where $\alpha, \beta, \sigma$ and $\omega$ satisfy the condition in Theorem 2. That generalize the initial data (\ref{A22}) stated in the introduction in some sense.
\item Condition appearing in the statement of Theorem 1 should be understood as follows: the quantity $\|e^{aD_{3}}u_{2}^{0}\|_{\dot{B}_{2,\gamma}^{\frac{1}{2\gamma^{\prime}},\frac{1}{2\gamma}}}$ and $\|u^{0}\|_{\dot{B}_{\infty,\infty}^{-1}}$ may be as large as wanted, as long as $\|e^{aD_{3}}u^{0}\|_{\dot{B}_{2,\gamma}^{\frac{1}{2\gamma^{\prime}}-\frac{1}{\gamma},\frac{1}{2\gamma}}}+\|e^{aD_{3}}u_{1}^{0}\|_{\dot{B}_{2,\gamma}^{\frac{1}{2\gamma^{\prime}},\frac{1}{2\gamma}}}
+\|e^{aD_{3}}\partial_{3}u_{3}^{0}\|_{\dot{B}_{2,\gamma}^{\frac{1}{2\gamma^{\prime}}-1,\frac{1}{2\gamma}}}$ is small enough. By using the assumption, one can also construct some other examples of large initial data. The problem we have to solve is the construction of large initial data $u^{0}$ such that a quadratic functional, namely (\ref{hang3332}) is small.
\end{itemize}
\begin{bf}Structure of the paper:\end{bf}
 In the second part, we recall some basic facts about Littlewood-Paley
 decomposition, anisotropic space and homogeneous Besov spaces, and we will also list some useful Propositions. The proof of Theorem 1 is achieved in section 6, using
some crucial Propositions. The proof of those Propositions is presented in Sections 4 and 5. This demands a careful use of the structure of nonlinear team for (NS). Theorem 2 is proved in subsection 6.3.\\\\
 \begin{bf}Notation.\end{bf} Throughout the paper, $\|\cdot\|_{L^{p}}$ denotes $L^{p}(\mathbb{R}^{d})$ norm of a function, $a\wedge b=min\{a,b\}$, $a\vee b=max\{a,b\}$ and $D_{i}f=(|\xi_{i}|\widehat{f})^{\vee}$. We shall sometimes use $X\lesssim Y$ to denote the
 estimate $X\leqslant CY$ for some C. In this paper, $C\geqslant1$, will denote constants which can be different at different places. The Fourier transform of $u$ is denoted either by $\widehat{u}$ or $\mathcal{F}u$, the inverse by $\mathcal{F}^{-1}u$.

 \section{Preliminaries}

Let us recall the Littlewood-Paley (or dyadic) decomposition. Let $\mathcal{S}$ be the Schwartz class of rapidly decreasing functions. Choose two nonnegative radial functions $\chi, \varphi\in\mathcal{S}(\mathbb{R}^{d})$, supported respectively in $\mathcal{B}=\{\xi\in \mathbb{R}^{d}: |\xi|\leqslant\frac{4}{3}\}$ and $\mathcal{C}=\{\xi\in \mathbb{R}^{d}: \frac{3}{4}\leqslant|\xi|\leqslant\frac{8}{3}\}$ such that
$$\chi(\xi)+\sum_{j\geqslant0}\varphi(2^{-j}\xi)=1,~\xi\in \mathbb{R}^{d},$$and
$$\sum_{j\in \mathbb{Z}}\varphi(2^{-j}\xi)=1,~ \xi\in \mathbb{R}^{d}\backslash \{0\}.$$
Denoting $\varphi_{j}(\xi)=\varphi(2^{-j}\xi)$ and $h_{j}=\mathcal{F}^{-1}\varphi_{j}$, we define the frequency localization operator as follows
$$\dot{\Delta}_{j}u:=\mathcal{F}^{-1}\varphi_{j}\mathcal{F}u=\int_{\mathbb{R}^{N}}h_{j}(y)u(x-y)dy,~~ \forall j\in \mathbb{Z},$$and
$$S_{j}u=\sum_{k\leqslant j-1}\dot{\Delta}_{k}u.$$
Informally, $\dot{\Delta}_{j}=S_{j}-S_{j-1}$ is a frequency projection to the annulus $\{|\xi|\thicksim2^{j}\}$, while $S_{j}$ is a frequency projection to the ball $\{|\xi|\lesssim2^{j}\}$. One easily obtains that
\begin{align}\label{LL1}
\dot{\Delta}_{j}\dot{\Delta}_{k}u\equiv0~~if~~|j-k|\geqslant2~~and~~\dot{\Delta}_{j}(S_{k-1}u\dot{\Delta}_{k}u)\equiv0~~if~~|j-k|\geqslant5.
\end{align}

Let us now recall the definition of homogeneous Besov spaces:\\\\
\begin{bf}Definition 2.1.\end{bf} For $s\in \mathbb{R}$, $1\leqslant p, r\leqslant\infty$ and $f\in \mathcal{Z}^{\prime}$, we set
$$\|f\|_{\dot{B}_{p,r}^{s}}:=(\sum_{j\in \mathbb{Z}}2^{jsr}\|\dot{\Delta}_{j}f\|_{L^{p}}^{r})^{\frac{1}{r}}<\infty,~~if~r<\infty,$$ and
 $$\|f\|_{\dot{B}_{p,\infty}^{s}}:=sup_{j\in \mathbb{Z}}2^{js}\|\dot{\Delta}_{j}f\|_{L^{p}}<\infty,$$
where $\mathcal{Z}^{\prime}$ denotes the dual space of $\mathcal{Z}=\{f\in \mathcal{S}(\mathbb{R}^{d}):\partial^{\gamma}\widehat{f}(0)=0;~\forall\gamma\in \mathbb{N}^{d}~multi-index\}$.

  The dyadic decomposition
operator is a frequency-localized operator, which is one of the most
delicate and clever ideas in harmonic analysis.
  As is well known, one of the interests of this decomposition is that the $\dot{\Delta}_{j}$ operators allow to count derivatives easily. More precisely, we recall the Bernstein inequality. A constant C exists such that $\forall k\in \mathbb{N}$, $1\leqslant p\leqslant q\leqslant\infty,$
 \begin{align}\label{hang114}sup_{|\alpha|=k}\|\partial^{\alpha}\dot{\Delta}_{j}u\|_{L^{q}(\mathbb{R}^{d})}\leqslant C^{k+1}2^{jk}2^{jd(\frac{1}{p}-\frac{1}{q})}\|\dot{\Delta}_{j}u\|_{L^{p}(\mathbb{R}^{d})}.\end{align}
\begin{bf}Proposition 2.2.\end{bf} The following properties on homogeneous Besov spaces hold true:

a. Embedding: Let $-\infty<s_{2}<s_{1}<\infty$, $s_{1}-\frac{d}{p}=s_{2}-\frac{d}{q}$ and $1\leqslant r, p, q\leqslant\infty$. Then $$\dot{B}_{p,r}^{s_{1}}(\mathbb{R}^{d})\hookrightarrow \dot{B}_{q,r}^{s_{2}}(\mathbb{R}^{d}).$$

b. Generalized derivative: Let $\sigma\in \mathbb{R}$, then the operator $D^{\sigma}$ is an isomorphism from $\dot{B}_{p,r}^{s}$ to $\dot{B}_{p,r}^{s-\sigma}$.

c. A constant c exists such that, for any integer $j$, any positive real number $t$ and any $p$ in $[1,\infty]$,
$$\|\dot{\Delta}_{j}e^{t\Delta}u\|_{L^{p}}\leqslant\frac{1}{c}e^{-c2^{2j}t}\|\dot{\Delta}_{j}u\|_{L^{p}}.$$

d. Equivalent definition for negative index Besov spaces: Let s be a positive real number, $1\leqslant p, r\leqslant\infty$. Denote
$$\{u\in \mathcal{Z}^{\prime}:\|u\|_{\dot{B}_{p,r}^{-s}}:=\|t^{\frac{s}{2}}\|e^{t\Delta}u\|_{L^{p}}\|_{L^{r}(\mathbb{R}^{+},\frac{dt}{t})}<\infty\}.$$
Proof. Using the Bernstein inequality (\ref{hang114}), it is easy to see that the first continuous embedding
holds. For the result of (b), one can refer to \cite{Triebel}, for the result of (c) and (d), one can refer to \cite{Chemin1}.

Below, let us recall the definition of nonhomogeneous Besov spaces:\\\\
\begin{bf}Definition 2.3.\end{bf} For $s\in \mathbb{R}$, $1\leqslant p, r\leqslant\infty$ and $f\in \mathcal{S}^{\prime}$, we set
$$\|f\|_{B_{p,r}^{s}}:=(\|\Delta_{0}f\|_{L^{p}}^{r}+\sum_{j>0}^{\infty}2^{jsr}\|\dot{\Delta}_{j}f\|_{L^{p}}^{r})^{\frac{1}{r}}<\infty.$$
where $\Delta_{0}f=\mathcal{F}^{-1}\chi\mathcal{F}f$ and $\mathcal{S}^{\prime}$ denotes the dual space of $\mathcal{S}$, which is said to be a tempered distribution space.
Note that we need replace
the $\ell^{r}$-norm by the $\ell^{\infty}$-norm in the above definition if $r=\infty$.

The frequency localization operators $\dot{\Delta}_{j}^{\nu}$ in the vertical direction are defined by
$$\dot{\Delta}_{j}^{\nu}u=\mathcal{F}_{\xi_{3}}^{-1}(\varphi(2^{-j}|\xi_{3}|)\mathcal{F}_{x_{i}}u)
:=\mathcal{F}_{\xi_{3}}^{-1}(\varphi_{j}^{\nu}\mathcal{F}_{x_{3}}u),~~S_{j}^{\nu}u=\sum_{k\leqslant j-1}\dot{\Delta}_{k}^{\nu}u,~~for~~j\in \mathbb{Z}.$$
Where $\mathcal{F}_{x_{3}}u=\int_{\mathbb{R}}u(x)e^{-2\pi ix_{3}\xi_{3}}dx_{3}$.
The frequency localization operators $\dot{\Delta}_{j}^{h}$ in the  horizontal direction are defined by
$$\dot{\Delta}_{j}^{h}u=\mathcal{F}_{\xi_{h}}^{-1}(\varphi(2^{-j}|\xi_{h}|)\mathcal{F}_{x_{h}}u)
:=\mathcal{F}_{\xi_{h}}^{-1}(\varphi_{j}^{h}\mathcal{F}_{x_{h}}u),~~S_{j}^{h}u=\sum_{k\leqslant j-1}\dot{\Delta}_{k}^{h}u,~~for~~j\in \mathbb{Z}.$$
Where $\mathcal{F}_{x_{h}}u=\int_{\mathbb{R}^{2}}u(x)e^{-2\pi i<x_{h},\xi_{h}>}dx_{h}$.

Now let us introduce the anisotropic homogeneous Besov spaces, which are important to control the low
horizontal frequency of the function. The role of anisotropic Sobolev or Besov
spaces in the study of a Navier-Stokes system with anisotropic vertical viscosity appears in \cite{Chemin14}.\\\\
\begin{bf}Definition 2.4.\end{bf} For $s_{1}, s_{2}\in \mathbb{R}$, $1\leqslant p, r\leqslant\infty$ and $f\in \mathcal{Z}^{\prime}$, we set
$$\|f\|_{\dot{B}_{p,r}^{s_{1},s_{2}}}:=(\sum_{i\in \mathbb{Z}}2^{is_{1}r}\sup_{j\in \mathbb{Z}}2^{js_{2}r}\|\dot{\Delta}_{i}^{h}\dot{\Delta}_{j}^{\nu}f\|_{L^{p}}^{r})^{\frac{1}{r}}<\infty,~~if~r<\infty,$$ and
 $$\|f\|_{\dot{B}_{p,\infty}^{s_{1},s_{2}}}:=sup_{i,j\in \mathbb{Z}}2^{is_{1}}2^{js_{2}}\|\dot{\Delta}_{i}^{h}\dot{\Delta}_{j}^{\nu}f\|_{L^{p}}<\infty,$$
where $\mathcal{Z}^{\prime}$ denotes the same notation as in Definition 2.1.

From the above definitions, one can easily see that $B_{p,r}^{s_{1}+s_{2}}\hookrightarrow\dot{B}_{p,r}^{s_{1}+s_{2}}\hookrightarrow\dot{B}_{p,r}^{s_{1},s_{2}}$, for $s_{1}, s_{2}>0$.
For $T>0$ and $\rho\in[1,\infty]$, we denote by $L^{\rho}([0,T),\dot{B}_{p,r}^{s_{1},s_{2}})$ the set of all tempered distribution $u$ satisfying
 $$\|u\|_{L^{\rho}([0,T),\dot{B}_{p,r}^{s_{1},s_{2}})}:=\|(\sum_{i\in \mathbb{Z}}2^{is_{1}r}\sup_{j\in \mathbb{Z}}2^{js_{2}r}\|\dot{\Delta}_{i}^{h}\dot{\Delta}_{j}^{\nu}u\|_{L^{p}}^{r})^{\frac{1}{r}}\|_{L_{T}^{\rho}}<\infty.$$
 We also need to use Chemin-Lerner-type Besov space $\widetilde{L}^{\rho}([0,T),\dot{B}_{p,r}^{s_{1},s_{2}})$, whose norm is defined by
 $$\|u\|_{\widetilde{L}^{\rho}([0,T),\dot{B}_{p,r}^{s_{1},s_{2}})}:=(\sum_{i\in \mathbb{Z}}2^{is_{1}r}\|\sup_{j\in \mathbb{Z}}2^{js_{2}r}\|\dot{\Delta}_{i}^{h}\dot{\Delta}_{j}^{\nu}u\|_{L^{p}}\|_{L_{T}^{\rho}}^{r})^{\frac{1}{r}}.$$
 For simplicity of notations, denote $L_{T}^{\rho}(\dot{B}_{p,r}^{s_{1},s_{2}}):=L^{\rho}([0,T),\dot{B}_{p,r}^{s_{1},s_{2}})$ and $\widetilde{L}_{T}^{\rho}(\dot{B}_{p,r}^{s_{1},s_{2}}):=\widetilde{L}^{\rho}([0,T),\dot{B}_{p,r}^{s_{1},s_{2}})$.
 By virtue of the Minkowski inequality, one has that
 $$\|u\|_{\widetilde{L}_{T}^{\rho}(\dot{B}_{p,r}^{s_{1},s_{2}})}\leqslant\|u\|_{L_{T}^{\rho}(\dot{B}_{p,r}^{s_{1},s_{2}})}~~if~\rho\leqslant r,$$
 and
 $$\|u\|_{L_{T}^{\rho}(\dot{B}_{p,r}^{s_{1},s_{2}})}\leqslant\|u\|_{\widetilde{L}_{T}^{\rho}(\dot{B}_{p,r}^{s_{1},s_{2}})}~~if~r\leqslant\rho.$$

Moreover, the proof of main results relies on exponential decay estimates for the Fourier transform of the solution. Thus, for any $\lambda\geqslant0$ and locally bounded function $\Phi$ on $\mathbb{R}^{+}\times\mathbb{R}^{3}$, and $\Psi$ on $\mathbb{R}^{+}$, and
for any function $f$, continuous in time and compactly supported in Fourier space, we define
$$f_{\lambda}(t):=\mathcal{F}^{-1}(e^{-\lambda\Psi(t)+\Phi(t,\cdot)}\widehat{f}(t,\cdot))~~~and~~~f_{\Phi}(t):=\mathcal{F}^{-1}(e^{\Phi(t,\cdot)}\widehat{f}(t,\cdot)).$$
Now we introduce two key quantities, which capture the weak
damping mechanism of the parabolic operator $\partial_{t}-\Delta$. We define the function $\theta(t)$ by
$$\dot{\theta}(t)=\|u_{\Phi}\|_{\dot{B}_{2,\gamma}^{\frac{1}{\gamma}+\frac{1}{2\gamma^{\prime}},\frac{1}{2\gamma}}}^{\gamma}~~and~~~\theta(0)=0,$$
and we also define
\begin{align}\label{12}\Psi(t)=\int_{0}^{t}\|(\bar{u}_{2})_{\Phi}(\tau)\|_{\dot{B}_{2,\gamma}^{\frac{5}{2\gamma^{\prime}},\frac{1}{2\gamma}}}^{\gamma^{\prime}}d\tau,\end{align}
where $\bar{u}_{2}$ satisfies Eq.(\ref{hang111}) and
\begin{align}\label{11}\Phi(t,\xi)=(a-k\theta(t))|\xi_{3}|\end{align}
 for some $k$ that will be chosen later on.

\section{The action of subadditive phases on products}

For any function $f$, we
denote by $f^{+}$ the inverse Fourier transform of $|\widehat{f}|$. Let us notice
that the map $f\mapsto f^{+}$ preserves the norm of all $\dot{B}_{2,q}^{s_{1},s_{2}}$ spaces.  Through out this section,
$\Phi$ will denote a locally bounded function on
$\mathbb{R}^{+}\times\mathbb{R}^{3}$ which satisfies
the following inequality: $$\Phi(t,\xi)\leqslant\Phi(t,\xi-\eta)+\Phi(t,\eta).$$

In what follows, we recall the paradifferential calculus which enables us to define a generalized product
between distributions, which is continuous in many functional spaces where the usual product
does not make sense.  The paraproduct between $f$ and $g$ is defined by
$T_{f}^{h}g:=\sum_{j\in\mathbb{Z}}S_{j-1}^{h}f\dot{\Delta}_{j}^{h}g$. We then have the following formal decomposition:
\begin{align}\label{14}fg=T_{f}^{h}g+T_{g}^{h}f+R^{h}(f,g)\end{align}
with$$R^{h}(f,g)=\sum_{j\in\mathbb{Z}}\dot{\Delta}_{j}^{h}f\tilde{\dot{\Delta}}_{j}^{h}g~~
and~~\tilde{\dot{\Delta}}^{h}_{j}=\dot{\Delta}^{h}_{j-1}+\dot{\Delta}^{h}_{j}+\dot{\Delta}^{h}_{j+1}.$$
The decomposition is called the Bony's paraproduct decomposition in the horizontal direction. In this paper, the decomposition will be used to study the nonlinear mapping estimates in anisotropic spaces.
Besides, we also need to use Bony's decomposition in
the vertical direction.\\\\
\begin{bf}Lemma 3.1.\end{bf} (Bernstein's inequality). Let $1\leqslant p\leqslant q\leqslant\infty$.  Assume that $f\in L^{p}(\mathbb{R}^{d})$. Then
there exists a constant C, independent of $f$, $j$, such that
$$supp\widehat{f}\subset \{|\xi|\leqslant C2^{j}\}\Rightarrow \|\partial^{\alpha}f\|_{L^{q}}\leqslant C2^{j|\alpha|+dj(\frac{1}{p}-\frac{1}{q})}\|f\|_{L^{p}},$$
$$supp\widehat{f}\subset \{\frac{1}{C}2^{j}\leqslant|\xi|\leqslant C2^{j}\}\Rightarrow C2^{-j|\alpha|}\sup_{|\beta|=|\alpha|}\|\partial^{\beta}f\|_{L^{p}}.$$
 \begin{bf}Lemma 3.2.\end{bf} Let $0<s_{1}<1$ and  $0<s_{2}<\frac{1}{2}$. And let $\gamma, \gamma_{1}, \gamma_{2}\in[1,\infty]$ with $\frac{1}{\gamma}=\frac{1}{\gamma_{1}}+\frac{1}{\gamma_{2}}$.  Assume that
$a_{\lambda}, b_{\lambda}\in \widetilde{L}_{T}^{\gamma_{1}}(\dot{B}_{2,1}^{s_{1},s_{2}})\cap \widetilde{L}_{T}^{\gamma_{2}}(\dot{B}_{2,1}^{1,\frac{1}{2}})$.  Then there holds
\begin{align*}&\|(ab)_{\lambda}\|_{\widetilde{L}_{T}^{\gamma}(\dot{B}_{2,1}^{s_{1},s_{2}})}\\&\leqslant e^{\lambda\Psi(T)}\|a_{\lambda}\|_{\widetilde{L}_{T}^{\gamma_{1}}(\dot{B}_{2,1}^{s_{1},s_{2}})}\|b_{\lambda}\|_{\widetilde{L}_{T}^{\gamma_{2}}(\dot{B}_{2,1}^{1,\frac{1}{2}})}
+e^{\lambda\Psi(T)}\|a_{\lambda}\|_{\widetilde{L}_{T}^{\gamma_{2}}(\dot{B}_{2,1}^{1,\frac{1}{2}})}\|b_{\lambda}\|_{\widetilde{L}_{T}^{\gamma_{1}}(\dot{B}_{2,1}^{s_{1},s_{2}})}.\end{align*}
\begin{bf}Proof\end{bf}.
Using Bony's decomposition (\ref{14}) in the horizontal direction, we write
\begin{align*}
\dot{\Delta}_{i}^{h}\dot{\Delta}_{j}^{\nu}(T_{a}^{\nu}b)&=\sum_{j^{\prime}}\dot{\Delta}_{i}^{h}\dot{\Delta}_{j}^{\nu}(S_{j^{\prime}-1}^{\nu}a\dot{\Delta}_{j^{\prime}}^{\nu}b)\\
&=\sum_{j^{\prime}}\dot{\Delta}_{i}^{h}\dot{\Delta}_{j}^{\nu}(T_{S_{j^{\prime}-1}^{\nu}a}^{h}\dot{\Delta}_{j^{\prime}}^{\nu}b+T_{\dot{\Delta}_{j^{\prime}}^{\nu}b}^{h}S_{j^{\prime}-1}^{\nu}a+R^{h}(S_{j^{\prime}-1}^{\nu}a,\dot{\Delta}_{j^{\prime}}^{\nu}b))\\
&:=A+B+C
\end{align*}
 Considering the support of the Fourier transform
of $T_{S_{j^{\prime}-1}^{\nu}a}^{h}\dot{\Delta}_{j^{\prime}}^{\nu}b$, we find
$$A=\sum_{|j-j^{\prime}|\leqslant4,|i-i^{\prime}|\leqslant4}\dot{\Delta}_{i}^{h}\dot{\Delta}_{j}^{\nu}(S_{i^{\prime}-1}^{h}S_{j^{\prime}-1}^{\nu}a\dot{\Delta}_{i^{\prime}}^{h}\dot{\Delta}_{j^{\prime}}^{\nu}b)$$
Then we get that
\begin{align}\label{hang553}
\|2^{is_{1}}2^{js_{2}}A\|_{L^{2}}&\lesssim\|S_{i^{\prime}-1}^{h}S_{j^{\prime}-1}^{\nu}a\|_{L^{\infty}}\sum_{|j-j^{\prime}|\leqslant4,|i-i^{\prime}|\leqslant4}\|2^{is_{1}}2^{js_{2}}\dot{\Delta}_{i^{\prime}}^{h}\dot{\Delta}_{j^{\prime}}^{\nu}b\|_{L^{2}}\nonumber\\
&\lesssim\|a\|_{\dot{B}_{2,1}^{s_{1},s_{2}}}\sum_{|j-j^{\prime}|\leqslant4,|i-i^{\prime}|\leqslant4}\|2^{i}2^{\frac{1}{2}j}\dot{\Delta}_{i^{\prime}}^{h}\dot{\Delta}_{j^{\prime}}^{\nu}b\|_{L^{2}}.
\end{align}
From the support of the Fourier transform
of $T_{\dot{\Delta}_{j^{\prime}}^{\nu}b}^{h}S_{j^{\prime}-1}^{\nu}a$, we have
$$B=\sum_{|j-j^{\prime}|\leqslant4,|i-i^{\prime}|\leqslant4}\dot{\Delta}_{i}^{h}\dot{\Delta}_{j}^{\nu}(\dot{\Delta}_{i^{\prime}}^{h}S_{j^{\prime}-1}^{\nu}aS_{i^{\prime}-1}^{h}\dot{\Delta}_{j^{\prime}}^{\nu}b).$$
 Then we get similarly that
 \begin{align}\label{hang554}
\|2^{is_{1}}2^{js_{2}}B\|_{L^{2}}&\lesssim\sum_{|j-j^{\prime}|\leqslant4\atop |i-i^{\prime}|\leqslant4}\|2^{is_{1}}\dot{\Delta}_{i^{\prime}}^{h}S_{j^{\prime}-1}^{\nu}a\|_{L_{x_{h}}^{2}L_{x_{3}}^{\infty}}\sum_{|j-j^{\prime}|\leqslant4\atop|i-i^{\prime}|\leqslant4}\|2^{js_{2}}S_{i^{\prime}-1}^{h}\dot{\Delta}_{j^{\prime}}^{\nu}b\|_{L_{x_{h}}^{\infty}L_{x_{3}}^{2}}\nonumber\\
&\lesssim\sum_{|i-i^{\prime}|\leqslant4}\|2^{is_{1}}2^{js_{2}}\dot{\Delta}_{i^{\prime}}^{h}\dot{\Delta}_{j}^{\nu}a\|_{\ell_{j}^{\infty}L^{2}}\sum_{|j-j^{\prime}|\leqslant4}\|2^{i}2^{\frac{1}{2}j}\dot{\Delta}_{i}^{h}\dot{\Delta}_{j^{\prime}}^{\nu}b\|_{\ell_{i}^{1}L^{2}}.
\end{align}
Now let us turn to C. We have
$$C=\sum_{|j-j^{\prime}|\leqslant4,i^{\prime}\geqslant i-2}\dot{\Delta}_{i}^{h}\dot{\Delta}_{j}^{\nu}(\dot{\Delta}_{i^{\prime}}^{h}S_{j^{\prime}-1}^{\nu}a\tilde{\dot{\Delta}}_{i^{\prime}}^{h}\dot{\Delta}_{j^{\prime}}^{\nu}b)$$
Then we get that
 \begin{align}\label{hang555}
&\|2^{is_{1}}2^{js_{2}}C\|_{L^{2}}\nonumber\\
&\lesssim\sum_{|j-j^{\prime}|\leqslant4,i^{\prime}\geqslant i-2}2^{(i-i^{\prime})s_{1}}2^{i^{\prime}s_{1}}\|\dot{\Delta}_{i^{\prime}}^{h}S_{j^{\prime}-1}^{\nu}a\|_{L^{\infty}}2^{js_{2}}\|\tilde{\dot{\Delta}}_{i^{\prime}}^{h}\dot{\Delta}_{j^{\prime}}^{\nu}b\|_{L^{2}}\\
&\lesssim\sum_{i^{\prime}\geqslant i-2}2^{(i-i^{\prime})s_{1}}2^{i^{\prime}s_{1}}\|2^{js_{2}}\dot{\Delta}_{i^{\prime}}^{h}\dot{\Delta}_{j}^{\nu}a\|_{\ell_{j}^{\infty}L^{2}}\sum_{|j-j^{\prime}|\leqslant4}2^{i^{\prime}}2^{\frac{1}{2}j}\|\dot{\Delta}_{i^{\prime}}^{h}\dot{\Delta}_{j^{\prime}}^{\nu}b\|_{L^{2}}.\nonumber
\end{align}
 Summing up (\ref{hang553})-(\ref{hang555}) yields that
 $$\|\sup_{j\in\mathbb{Z}}2^{is_{1}}2^{js_{2}}\|\dot{\Delta}_{i}^{h}\dot{\Delta}_{j}^{\nu}(T_{a}^{\nu}b)\|_{L^{2}}\|_{\ell_{i}^{1}L_{T}^{\gamma}}\lesssim\|a\|_{\widetilde{L}_{T}^{\gamma_{1}}(\dot{B}_{2,1}^{s_{1},s_{2}})}\|b\|_{\widetilde{L}_{T}^{\gamma_{2}}(\dot{B}_{2,1}^{1,\frac{1}{2}})}.$$
From a very similar arguments as in the above, one has that
$$\|\sup_{j\in\mathbb{Z}}2^{is_{1}}2^{js_{2}}\|\dot{\Delta}_{i}^{h}\dot{\Delta}_{j}^{\nu}[(T_{b}^{\nu}a)+R^{\nu}(a,b)]\|_{L^{2}}\|_{\ell_{i}^{1}L_{T}^{\gamma}}\lesssim
\|a\|_{\widetilde{L}_{T}^{\gamma_{2}}(\dot{B}_{2,1}^{1,\frac{1}{2}})}\|b\|_{\widetilde{L}_{T}^{\gamma_{1}}(\dot{B}_{2,1}^{s_{1},s_{2}})}.$$
The lemma is proved in the case when the
function $-\lambda\Psi+\Phi$ is identically 0. In order to deal with the general case, we only need
to notice the fact that
$$|\mathcal{F}[\dot{\Delta}_{i}^{h}\dot{\Delta}_{j}^{\nu}ab]_{\lambda}(\xi)|\leqslant e^{\lambda\Psi(T)}|\mathcal{F}[\dot{\Delta}_{i}^{h}\dot{\Delta}_{j}^{\nu}a_{\lambda}^{+}b_{\lambda}^{+}](\xi)|.$$
The proof of Lemma 3.2 is finished.

From a similar arguments, one can also obtain the estimates for $(ab)_{\lambda}$ in $\widetilde{L}_{T}^{\gamma}(\dot{B}_{2,q}^{s_{1},s_{2}})$ with $q>1$.

\section{The action of the phase $-\lambda\Psi+\Phi$ on the heat operator}
Applying the Duhamel formula to Eq.(\ref{hang113}) gives
\begin{align}\label{21}u(t)=e^{t\Delta}u^{0}-\int_{0}^{t}e^{(t-\tau)\Delta}\mathbb{P}[(u\cdot\nabla)u]d\tau,\end{align}
 where $e^{t\Delta}=\mathcal{F}^{-1}e^{-t|\xi|^{2}}\mathcal{F}$ and $$[\mathbb{P}((u\cdot\nabla)u)]_{i}=(u\cdot\nabla)u_{i}+R_{i}\sum_{k=1}^{3}R_{k}[(u\cdot\nabla)u_{k}].$$
 For simplicity, denote
 $$(\mathcal{A}f)(t,x):=\int_{0}^{t}e^{(t-\tau)\Delta}f(\tau,x)d\tau.$$

In $e^{t\Delta}$, $e^{-t|\xi|^{2}}$ as an exponential
 decay function, which corresponds to the dissipation for
 $e^{t\Delta}$, determines that the semi-group $e^{t\Delta}$
 has very good properties. The
heat semi-group has been widely studied in the mixes space
$L^{\gamma}(\mathbb{R}^{+};\dot{B}_{p,q}^{s})$ and the technique is to
extensively use the dyadic decomposition together with the
exponential decay of $e^{-t|\xi|^{2}}$, see \cite{chen1}, \cite{Wang}, \cite{Ru}. Let $\Psi$ and $\Phi(t,\xi)$ be given by
(\ref{11}) and (\ref{12}). In this section, we will study the action
of the multiplier $e^{-\lambda\Psi(t)+\Phi(t,\xi)}$ on the heat operator $e^{t\Delta}$.\\\\
\begin{bf}Lemma 4.1.\end{bf} Let $s_{1}, s_{2}\geqslant0$, $s,\sigma\in\mathbb{R}$, $2\leqslant p\leqslant\infty$, $1\leqslant q\leqslant\infty$. Then we have
$$\|(e^{t\Delta}u^{0})_{\lambda}\|_{\widetilde{L}_{T}^{\gamma}(\dot{B}_{p,q}^{s,\sigma}(\mathbb{R}^{3}))}\lesssim\|e^{aD_{3}}u^{0}\|_{\dot{B}_{2,q}^{s+1-\frac{2}{p}-\frac{2}{\gamma},\sigma+\frac{1}{2}-\frac{1}{p}}(\mathbb{R}^{3})}.$$
Proof. Thanks to (C) in Pro.2.2 and the definition of $\widetilde{L}_{T}^{\gamma}(\dot{B}_{p,q}^{s,\sigma}(\mathbb{R}^{3}))$, we have
\begin{align*}
&\|e^{-\lambda\Psi(t)+\Phi(t,D)}e^{t\Delta}u^{0}\|_{\widetilde{L}_{T}^{\gamma}(\dot{B}_{p,q}^{s,\sigma})}\\
&\lesssim
\|2^{is}\sup_{j\in\mathbb{Z}}2^{j\sigma}\|\dot{\Delta}_{i}^{h}\dot{\Delta}_{j}^{\nu}e^{-\lambda\Psi(t)+\Phi(t,D)}e^{-ct(2^{2i}+2^{2j})}u^{0}\|_{L^{p}}\|_{\ell_{i}^{q}L_{T}^{\gamma}}\\
&\lesssim(\sum_{i\in\mathbb{Z}}2^{iq(s+1-\frac{2}{p})}\|e^{-ct2^{2i}}\|_{L_{T}^{\gamma}}^{q}\sup_{j\in\mathbb{Z}}2^{jq(\sigma+\frac{1}{2}-\frac{1}{p})}\|\dot{\Delta}_{i}^{h}\dot{\Delta}_{j}^{\nu}e^{aD_{3}}u^{0}\|_{L^{2}}^{q})^{\frac{1}{q}}\\
&\lesssim(\sum_{i\in\mathbb{Z}}2^{i(s+1-\frac{2}{p}-\frac{2}{\gamma})q}\sup_{j\in\mathbb{Z}}2^{jq(\sigma+\frac{1}{2}-\frac{1}{p})}\|\dot{\Delta}_{i}^{h}\dot{\Delta}_{j}^{\nu}e^{aD_{3}}u^{0}\|_{L^{2}}^{q})^{\frac{1}{q}}\\
&\leqslant\|e^{aD_{3}}u^{0}\|_{\dot{B}_{2,q}^{s+1-\frac{2}{p}-\frac{2}{\gamma},\sigma+\frac{1}{2}-\frac{1}{p}}}.
\end{align*}
The proof of Lemma 4.1 is finished.\\\\
\begin{bf}Lemma 4.2.\end{bf} Let $s_{1}, s_{2}, s_{1}^{\prime}, s_{2}^{\prime}\in \mathbb{R}$ and $1\leqslant \gamma,
\gamma_{1}^{\prime},q\leqslant
\infty$. Assume that $\gamma_{1}^{\prime}\leqslant\gamma$, $s_{1}^{\prime}<s_{1}$ and $s_{1}^{\prime}+s_{2}^{\prime}=s_{1}+s_{2}-\frac{2}{\gamma}-\frac{2}{\gamma_{1}}$. Then, for any $T>0$, we have
$$\|e^{-\lambda\Psi(t)+\Phi(t,D)}\mathcal{A}f\|_{\widetilde{L}_{T}^{\gamma}(\dot{B}_{2,q}^{s_{1},s_{2}}(\mathbb{R}^{3}))}
\lesssim\|f_{\lambda}\|_{\widetilde{L}_{T}^{\gamma_{1}^{\prime}}(\dot{B}_{2,q}^{s_{1}^{\prime},s_{2}^{\prime}}(\mathbb{R}^{3}))}.$$
Proof. Applying the Fourier multiplier $e^{-\lambda\Psi(t)+\Phi(t,D)}\dot{\Delta}_{i}^{h}\dot{\Delta}_{j}^{\nu}$ to $\mathcal{A}f$ gives
\begin{align*}
e^{-\lambda\Psi(t)+\Phi(t,D)}\dot{\Delta}_{i}^{h}\dot{\Delta}_{j}^{\nu}\mathcal{A}f=\int_{0}^{t}e^{(t-\tau)\Delta-\lambda\Psi(t)+\Phi(t,D)+\lambda\Psi(\tau)-\Phi(\tau,D)}(\dot{\Delta}_{i}^{h}\dot{\Delta}_{j}^{\nu}f)_{\lambda}d\tau.
\end{align*}
 Due to the definition of $\Phi$ and $\Psi$, we have
 \begin{align*}&-\lambda\Psi(t)+\Phi(t,\xi)+\lambda\Psi(\tau)-\Phi(\tau,\xi)\\&=-\lambda\int_{\tau}^{t}\|(\bar{u}_{2})_{\Phi}(\tau)\|_{\dot{B}_{2,\gamma}^{\frac{5}{2\gamma^{\prime}},\frac{1}{2\gamma}}}^{\gamma^{\prime}}d\tau-k|\xi_{3}|\int_{\tau}^{t}\dot{\theta}(\tau^{\prime})d\tau^{\prime}.
 \end{align*}
 Then we take the Fourier transform to get
\begin{align*}
\mathcal{F}e^{-\lambda\Psi(t)+\Phi(t,D)}\dot{\Delta}_{i}^{h}\dot{\Delta}_{j}^{\nu}\mathcal{A}f=\int_{0}^{t}e^{-(t-\tau)|\xi|^{2}-\lambda\Psi(t)+\Phi(t,\xi)+\lambda\Psi(\tau)-\Phi(\tau,\xi)}\mathcal{F}(\dot{\Delta}_{i}^{h}\dot{\Delta}_{j}^{\nu}f)_{\lambda}d\tau.
\end{align*}
We take the $L^{2}$ norm to obtain
\begin{align}\label{22}
\|e^{-\lambda\Psi(t)+\Phi(t,D)}\dot{\Delta}_{i}^{h}\dot{\Delta}_{j}^{\nu}\mathcal{A}f\|_{L^{2}}\leqslant\int_{0}^{t}\|e^{-c(t-\tau)(2^{2i}+2^{2j})}(\dot{\Delta}_{i}^{h}\dot{\Delta}_{j}^{\nu}f)_{\lambda}\|_{L^{2}}d\tau.
\end{align}
Noticing that
$$2^{j(s_{2}-s_{2}^{\prime})}e^{-t2^{2j}}\leqslant t^{-\frac{s_{2}-s_{2}^{\prime}}{2}}.$$
From which and (\ref{22}) we infer that
\begin{align}\label{23}
&\|\sup_{j\in\mathbb{Z}}2^{s_{2}j}\|e^{-\lambda\Psi(t)+\Phi(t,D)}\dot{\Delta}_{i}^{h}\dot{\Delta}_{j}^{\nu}\mathcal{A}f\|_{L^{2}}\|_{L_{T}^{\gamma}}\nonumber\\
&\lesssim\|\int_{0}^{t}e^{-c(t-\tau)2^{2i}}(t-\tau)^{-\frac{s_{2}-s_{2}^{\prime}}{2}}\sup_{j\in\mathbb{Z}}2^{s_{2}^{\prime}j}\|(\dot{\Delta}_{i}^{h}\dot{\Delta}_{j}^{\nu}f)_{\lambda}\|_{L^{2}}d\tau\|_{L_{T}^{\gamma}}.
\end{align}
Due to $s_{2}-s_{2}^{\prime}<\frac{2}{\gamma}+\frac{2}{\gamma_{1}}$, we have
$$\int_{0}^{t}e^{-c\tau2^{2i}}\tau^{-\frac{(s_{2}-s_{2}^{\prime})\gamma\gamma_{1}}{2(\gamma+\gamma_{1})}}d\tau\lesssim 2^{2(-1+\frac{(s_{2}-s_{2}^{\prime})\gamma\gamma_{1}}{2(\gamma+\gamma_{1})})i}.$$
 Then, applying Young's inequality gives
\begin{align*}
&(\sum_{i\in\mathbb{Z}}2^{iqs_{1}}\|\sup_{j\in\mathbb{Z}}2^{s_{2}j}\|e^{-\lambda\Psi(t)+\Phi(t,D)}\dot{\Delta}_{i}^{h}\dot{\Delta}_{j}^{\nu}\mathcal{A}f\|_{L^{2}}\|_{L_{T}^{\gamma}}^{q})^{\frac{1}{q}}\\
&\lesssim(\sum_{i\in\mathbb{Z}}2^{s_{1}^{\prime}iq}\|\sup_{j\in\mathbb{Z}}2^{s_{2}^{\prime}j}\|(\dot{\Delta}_{i}^{h}\dot{\Delta}_{j}^{\nu}f)_{\lambda}\|_{L^{2}}\|_{L_{T}^{\gamma_{1}^{\prime}}}^{q})^{\frac{1}{q}}.
\end{align*}
 This gives the inequality of Lemma 4.2.

Based on the above Lemmas 4.1 and 4.2,
we can obtain the following A priori estimates for the Eq.(\ref{hang111}). For simplicity, we denote $X_{T}:=\widetilde{L}_{T}^{\infty}(\dot{B}_{2,\gamma}^{\frac{1}{2\gamma^{\prime}},\frac{1}{2\gamma}})\cap \widetilde{L}_{T}^{\gamma^{\prime}}(\dot{B}_{2,\gamma}^{\frac{5}{2\gamma^{\prime}},\frac{1}{2\gamma}})$, $Y=\dot{B}_{2,\gamma}^{\frac{5}{2\gamma^{\prime}},\frac{1}{2\gamma}}$ and $Z=\dot{B}_{2,\gamma}^{\frac{1}{\gamma}+\frac{1}{2\gamma^{\prime}},\frac{1}{2\gamma}}$.\\\\
\begin{bf}Proposition 4.3.\end{bf} Let $\lambda, k>0$ and $1<\gamma<2$. A constant C exists such that for any $T>0$
satisfying $\theta(t)\leqslant a/k$ for $t\in[0,T]$, we have
\begin{align*}&\|\mathcal{A}div_{h}(\bar{u}_{h}\bar{u}_{1})_{\lambda}\|_{X_{T}}+\|\Lambda_{h}^{-1}\partial_{2}\mathcal{A}div_{h}(\bar{u}_{h}\bar{u}_{2})_{\lambda}\|_{X_{T}}+\|\Lambda_{h}^{-1}\partial_{3}\mathcal{A}div(uu_{3})_{\lambda}\|_{X_{T}}\\
&\leqslant \frac{C}{k^{\frac{1}{\gamma}}}
\|(u_{h})_{\lambda}\|_{X_{T}}+e^{\lambda\Psi(T)}
\|(\bar{u}_{1})_{\lambda}\|_{X_{T}}\|(\bar{u}_{1})_{\lambda}\|_{X_{T}}\\
&~~~+\frac{C}{\lambda^{^{\frac{1}{\gamma^{\prime}}}}}
(\|(\partial_{2}\bar{u}_{2})_{\lambda}\|_{\widetilde{L}_{T}^{\gamma}(\dot{B}_{2,\gamma}^{\frac{2}{\gamma}+\frac{1}{2\gamma^{\prime}}-1,\frac{1}{2\gamma}})}
+\|(\bar{u}_{1})_{\lambda}\|_{\widetilde{L}_{T}^{\gamma^{\prime}}(\dot{B}_{2,\gamma}^{\frac{5}{2\gamma^{\prime}},\frac{1}{2\gamma}})}).
\end{align*}
\begin{bf}Proof.\end{bf} Noticing that
$$2^{\frac{3}{2}j}e^{-t2^{2j}}\leqslant C t^{-\frac{3}{4}},$$
 we infer from (c) in Proposition 2.2 and Young's inequality that, for $1<\gamma<2$,
\begin{align*}
&\|\Lambda_{h}^{-1}\partial_{3}e^{-\lambda\Psi(t)+\Phi(t,D)}\mathcal{A}(
u_{3}div_{h}u_{h})\|_{\widetilde{L}_{T}^{\infty}(\dot{B}_{2,\gamma}^{\frac{1}{2\gamma^{\prime}},\frac{1}{2\gamma}})}\\
&\lesssim\|\int_{0}^{t}e^{-c(t-\tau)(2^{2i}+2^{2j})-k2^{j}\int_{\tau}^{t}\dot{\theta}(\tau^{\prime})d\tau^{\prime}}\\
&~~~~~~~\times
\frac{2^{\frac{1}{2\gamma^{\prime}} i}2^{(\frac{1}{2\gamma}+1)j}\|\dot{\Delta}_{i}^{h}\dot{\Delta}_{j}^{\nu}\Lambda_{h}^{-1}(u_{3}
div_{h}u_{h})_{\lambda}(\tau,\cdot)\|_{L^{2}}\|(u_{3})_{\Phi}\|_{\dot{B}_{2,\gamma}^{\frac{1}{\gamma}+\frac{1}{2\gamma^{\prime}},\frac{1}{2\gamma}}}}{\|(u_{3})_{\Phi}\|_{\dot{B}_{2,\gamma}^{\frac{1}{\gamma}+\frac{1}{2\gamma^{\prime}},\frac{1}{2\gamma}}}}d\tau\|_{\ell_{i}^{\gamma}L_{T}^{\infty}\ell_{j}^{\infty}}\nonumber\\
&\lesssim
\|(\int_{0}^{t}e^{-c\gamma^{\prime}(t-\tau)2^{2i}}(t-\tau)^{-\frac{3}{4}}\frac{2^{\frac{1}{2}i}2^{(\frac{1}{2\gamma}-\frac{1}{2\gamma^{\prime}})\gamma^{\prime} j}\|\dot{\Delta}_{i}^{h}\dot{\Delta}_{j}^{\nu}\Lambda_{h}^{-1}(u_{3}
div_{h}u_{h})_{\lambda}(\tau,\cdot)\|_{L^{2}}^{\gamma^{\prime}}}{\|(u_{3})_{\Phi}\|_{\dot{B}_{2,\gamma}^{\frac{1}{\gamma}+\frac{1}{2\gamma^{\prime}},\frac{1}{2\gamma}}}^{\gamma^{\prime}}}d\tau)^{\frac{1}{\gamma^{\prime}}}\nonumber\\
&~~~~~~\times(\int_{0}^{t}e^{-\gamma k2^{j}\int_{\tau}^{t}\dot{\theta}(\tau^{\prime})d\tau^{\prime}}
2^{j}\|(u_{3})_{\Phi}\|_{\dot{B}_{2,\gamma}^{\frac{1}{\gamma}+\frac{1}{2\gamma^{\prime}},\frac{1}{2\gamma}}}^{\gamma}
d\tau)^{\frac{1}{\gamma}}\|_{\ell_{i}^{\gamma}L_{T}^{\infty}\ell_{j}^{\infty}}.
\end{align*}
From the proof of Lemma 3.2, it is easy to see that
 \begin{align}\label{19}
  &2^{\frac{1}{2\gamma^{\prime}}i}2^{(\frac{1}{2\gamma}-\frac{1}{2\gamma^{\prime}})j}\|\dot{\Delta}_{i}^{h}\dot{\Delta}_{j}^{\nu}\Lambda_{h}^{-1}(u_{3}
div_{h}u_{h})_{\lambda}(\tau,\cdot)\|_{L^{2}}\nonumber\\
&\lesssim 2^{(-\frac{1}{\gamma}-1)i}\|(u_{3})_{\Phi}\|_{\dot{B}_{2,\gamma}^{\frac{1}{\gamma}+\frac{1}{2\gamma^{\prime}},\frac{1}{2\gamma}}}\sum_{i^{\prime}\leqslant i+4}2^{i^{\prime}}\sup_{j\in\mathbb{Z}}2^{\frac{1}{2\gamma}j}
\|\dot{\Delta}_{i^{\prime}}^{h}\dot{\Delta}_{j}^{\nu}(div_{h}u_{h})_{\lambda}\|_{L^{2}}\\
&~~~~+2^{\frac{1}{2\gamma^{\prime}}i}\|(u_{3})_{\Phi}\|_{\dot{B}_{2,\gamma}^{\frac{1}{\gamma}+\frac{1}{2\gamma^{\prime}},\frac{1}{2\gamma}}}\sum_{i^{\prime}\geqslant i-4}2^{\frac{1}{2\gamma^{\prime}}i^{\prime}}\sup_{j\in\mathbb{Z}}2^{\frac{1}{2\gamma}j}
\|\dot{\Delta}_{i^{\prime}}^{h}\dot{\Delta}_{j}^{\nu}(u_{h})_{\lambda}\|_{L^{2}}.\nonumber
 \end{align}
 Consequently, we infer from
 $$\int_{0}^{t}e^{-c\tau2^{2i}}\tau^{-\frac{3}{4}\cdot\frac{8}{7}}d\tau\lesssim 2^{-\frac{2}{7}i},$$
 that
\begin{align}\label{hang51}
&\|\Lambda_{h}^{-1}\partial_{3}e^{-\lambda\Psi(t)+\Phi(t,D)}\mathcal{A}(
u_{3}div_{h}u_{h})\|_{\widetilde{L}_{T}^{\infty}(\dot{B}_{2,\gamma}^{\frac{1}{2\gamma^{\prime}},\frac{1}{2\gamma}})}\nonumber\\
&\lesssim\frac{1}{k^{\frac{1}{\gamma}}}
\|(\int_{0}^{t}e^{-c\gamma^{\prime}\tau2^{2i}}\tau^{-\frac{3}{4}\cdot\frac{8}{7}}d\tau)^{\frac{7}{8\gamma^{\prime}}}\\
&~~~\times(\int_{0}^{T}\frac{2^{4i}\sup_{j\in\mathbb{Z}}2^{8(\frac{1}{2\gamma}-\frac{1}{2\gamma^{\prime}})\gamma^{\prime} j}\|\dot{\Delta}_{i}^{h}\dot{\Delta}_{j}^{\nu}\Lambda_{h}^{-1}(u_{3}
div_{h}u_{h})_{\lambda}(\tau,\cdot)\|_{L^{2}}^{8\gamma^{\prime}}}{\|(u_{3})_{\Phi}\|_{\dot{B}_{2,\gamma}^{\frac{1}{\gamma}+\frac{1}{2\gamma^{\prime}},\frac{1}{2\gamma}}}^{8\gamma^{\prime}}}d\tau)^{\frac{1}{8\gamma^{\prime}}}\|_{\ell_{i}^{\gamma}}\nonumber\\
&\lesssim\frac{1}{k^{\frac{1}{\gamma}}}
\|(u_{h})_{\lambda}\|_{\widetilde{L}_{T}^{8\gamma^{\prime}}(\dot{B}_{2,\gamma}^{\frac{3}{4\gamma^{\prime}},\frac{1}{2\gamma}})}.\nonumber
\end{align}
Similarly, we have
\begin{align}\label{hang5151}
&\|\Lambda_{h}^{-1}\partial_{3}e^{-\lambda\Psi(t)+\Phi(t,D)}\mathcal{A}div_{h}(
u_{h}u_{3})\|_{\widetilde{L}_{T}^{\infty}(\dot{B}_{2,\gamma}^{\frac{1}{2\gamma^{\prime}},\frac{1}{2\gamma}})}\nonumber\\
&\lesssim
\|(\int_{0}^{t}e^{-c\gamma^{\prime}(t-\tau)2^{2i}}(t-\tau)^{-\frac{3}{4}}\frac{2^{\frac{1}{2}i}2^{(\frac{1}{2\gamma}-\frac{1}{2\gamma^{\prime}})\gamma^{\prime} j}\|\dot{\Delta}_{i}^{h}\dot{\Delta}_{j}^{\nu}(u_{h}
u_{3})_{\lambda}(\tau,\cdot)\|_{L^{2}}^{\gamma^{\prime}}}{\|(u_{3})_{\Phi}\|_{Z}^{\gamma^{\prime}}}d\tau)^{\frac{1}{\gamma^{\prime}}}\\
&~~~~~~~~~~\times(\int_{0}^{t}e^{-\gamma k2^{j}\int_{\tau}^{t}\dot{\theta}(\tau^{\prime})d\tau^{\prime}}
2^{j}\|(u_{3})_{\Phi}\|_{Z}^{\gamma}
d\tau)^{\frac{1}{\gamma}}\|_{\ell_{i}^{\gamma}L_{T}^{\infty}\ell_{j}^{\infty}}\nonumber\\
&\lesssim\frac{1}{k^{\frac{1}{\gamma}}}
\|(u_{h})_{\lambda}\|_{\widetilde{L}_{T}^{8\gamma^{\prime}}(\dot{B}_{2,\gamma}^{\frac{3}{4\gamma^{\prime}},\frac{1}{2\gamma}})}.\nonumber
\end{align}
On the other hand, noticing that
$$2^{\frac{3}{2\gamma^{\prime}}j}e^{-t2^{2j}}\leqslant C t^{-\frac{3}{4\gamma^{\prime}}},$$
and
$$\int_{0}^{t}e^{-c\tau2^{2i}}\tau^{-\frac{3}{4}}d\tau\lesssim 2^{-\frac{1}{2}i},$$
similar to (\ref{hang51}), (\ref{hang5151}) and (\ref{23}), applying Young's inequality gives
\begin{align}\label{hang52}
&\|\Lambda_{h}^{-1}\partial_{3}e^{-\lambda\Psi(t)+\Phi(t,D)}\mathcal{A}(
u_{3}div_{h}u_{h})\|_{\widetilde{L}_{T}^{\gamma^{\prime}}(\dot{B}_{2,\gamma}^{\frac{5}{2\gamma^{\prime}},\frac{1}{2\gamma}})}\nonumber\\
&\lesssim\|(\int_{0}^{t}(t-\tau)^{-\frac{3}{4}}e^{-\gamma^{\prime}c(t-\tau)2^{2i}}\frac{2^{\frac{5}{2}i}2^{(\frac{2}{\gamma}-\frac{1}{2})\gamma^{\prime} j}2^{-\frac{\gamma^{\prime}}{\gamma} j}\|\dot{\Delta}_{i}^{h}\dot{\Delta}_{j}^{\nu}\Lambda_{h}^{-1}(u_{3}div_{h}u_{h})_{\lambda}(\tau,\cdot)\|_{L^{2}}^{\gamma^{\prime}}}{\|\bar{u}_{3}\|_{Z}^{\gamma^{\prime}}}
d\tau)^{\frac{1}{\gamma^{\prime}}}\nonumber\\
&~~~~~~~~~~~\times(\int_{0}^{t}e^{-\gamma k2^{j}\int_{\tau}^{t}\dot{\theta}(\tau^{\prime})d\tau^{\prime}}
2^{j}\|(u_{3})_{\Phi}\|_{Z}^{\gamma}
d\tau)^{\frac{1}{\gamma}}\|_{\ell_{i}^{\gamma}L_{T}^{\gamma^{\prime}}\ell_{j}^{\infty}}\\
&\lesssim\frac{1}{k^{\frac{1}{\gamma}}}
\|(u_{h})_{\lambda}\|_{\widetilde{L}_{T}^{\gamma^{\prime}}(\dot{B}_{2,\gamma}^{\frac{5}{2\gamma^{\prime}},\frac{1}{2\gamma}})}.\nonumber
\end{align}
and
\begin{align}\label{hang152}
&\|\Lambda_{h}^{-1}\partial_{3}e^{-\lambda\Psi(t)+\Phi(t,D)}\mathcal{A}div_{h}(
u_{h}u_{3})\|_{\widetilde{L}_{T}^{\gamma^{\prime}}(\dot{B}_{2,\gamma}^{\frac{5}{2\gamma^{\prime}},\frac{1}{2\gamma}})}\lesssim\frac{1}{k^{\frac{1}{\gamma}}}
\|(u_{h})_{\lambda}\|_{\widetilde{L}_{T}^{\gamma^{\prime}}(\dot{B}_{2,\gamma}^{\frac{5}{2\gamma^{\prime}},\frac{1}{2\gamma}})}.
\end{align}
 Due to $$2^{\frac{\gamma}{2\gamma^{\prime}}j}e^{-t2^{2j}}\leqslant C t^{-\frac{\gamma}{4\gamma^{\prime}}},$$ from a similar arguments as in the above, one can obtain that
\begin{align}\label{hang92}
&\|e^{-\lambda\Psi(t)+\Phi(t,D)}\mathcal{A}\partial_{2}(
\bar{u}_{2}
\bar{u}_{2})\|_{\widetilde{L}_{T}^{\gamma^{\prime}}(\dot{B}_{2,\gamma}^{\frac{5}{2\gamma^{\prime}},\frac{1}{2\gamma}})}\nonumber\\
&\lesssim
\|(\int_{0}^{t}(t-\tau)^{-\frac{\gamma}{4\gamma^{\prime}}}e^{-c\gamma(t-\tau)2^{2i}}\frac{2^{\frac{5\gamma}{2\gamma^{\prime}}i}2^{\frac{1}{2}(\frac{1}{\gamma}-\frac{1}{\gamma^{\prime}})\gamma j}\|\dot{\Delta}_{i}^{h}\dot{\Delta}_{j}^{\nu}\partial_{2}(
\bar{u}_{2}
\bar{u}_{2})_{\lambda}(\tau,\cdot)\|_{L^{2}}^{\gamma}}{\|(\bar{u}_{2})_{\Phi}\|_{\dot{B}_{2,\gamma}^{\frac{5}{2\gamma^{\prime}},\frac{1}{2\gamma}}}^{\gamma}}d\tau)^{\frac{1}{\gamma}}\nonumber\\
&~~~~~~~~~~\times(\int_{0}^{t}e^{-\gamma^{\prime} \lambda\int_{\tau}^{t}\|(\bar{u}_{2})_{\Phi}(\tau^{\prime})\|_{Y}^{\gamma^{\prime}}d\tau^{\prime}}
\|(\bar{u}_{2})_{\Phi}\|_{\dot{B}_{2,\gamma}^{\frac{5}{2\gamma^{\prime}},\frac{1}{2\gamma}}}^{\gamma^{\prime}}
d\tau)^{\frac{1}{\gamma^{\prime}}}\|_{\ell_{i}^{\gamma}L_{T}^{\gamma^{\prime}}\ell_{j}^{\infty}}\\
&\leqslant\frac{C}{\lambda^{^{\frac{1}{\gamma^{\prime}}}}}
\|(\partial_{2}\bar{u}_{2})_{\lambda}\|_{\widetilde{L}_{T}^{\gamma}(\dot{B}_{2,\gamma}^{\frac{2}{\gamma}+\frac{1}{2\gamma^{\prime}}-1,\frac{1}{2\gamma}})}.\nonumber
\end{align}
Similarly, one has that
\begin{align}\label{hang192}
&\|e^{-\lambda\Psi(t)+\Phi(t,D)}\mathcal{A}[\partial_{2}(
\bar{u}_{1}
\bar{u}_{2})+\partial_{1}(
\bar{u}_{1}\bar{u}_{1})]\|_{\widetilde{L}_{T}^{\gamma^{\prime}}(\dot{B}_{2,\gamma}^{\frac{5}{2\gamma^{\prime}},\frac{1}{2\gamma}})}\nonumber\\
&\lesssim\frac{C}{\lambda^{^{\frac{1}{\gamma^{\prime}}}}}
\|(\bar{u}_{1})_{\lambda}\|_{\widetilde{L}_{T}^{\gamma^{\prime}}(\dot{B}_{2,\gamma}^{\frac{5}{2\gamma^{\prime}},\frac{1}{2\gamma}})}
+e^{\lambda\Psi(T)}\|(\bar{u}_{1})_{\lambda}\|_{X_{T}}\|(\bar{u}_{1})_{\lambda}\|_{\widetilde{L}_{T}^{\gamma^{\prime}}(\dot{B}_{2,\gamma}^{\frac{5}{2\gamma^{\prime}},\frac{1}{2\gamma}})}
\end{align}
and
\begin{align}\label{hang1192}
&\|e^{-\lambda\Psi(t)+\Phi(t,D)}\mathcal{A}[\partial_{2}(
\bar{u}_{2}
\bar{u}_{2})+\partial_{2}(
\bar{u}_{1}
\bar{u}_{2})+\partial_{1}(
\bar{u}_{1}\bar{u}_{1})]\|_{\widetilde{L}_{T}^{\infty}(\dot{B}_{2,\gamma}^{\frac{1}{2\gamma^{\prime}},\frac{1}{2\gamma}})}\nonumber\\
&\leqslant\frac{C}{\lambda^{^{\frac{1}{\gamma^{\prime}}}}}
(\|(\partial_{2}\bar{u}_{2})_{\lambda}\|_{\widetilde{L}_{T}^{\gamma}(\dot{B}_{2,\gamma}^{\frac{2}{\gamma}+\frac{1}{2\gamma^{\prime}}-1,\frac{1}{2\gamma}})}+\|(\bar{u}_{1})_{\lambda}\|_{\widetilde{L}_{T}^{\gamma^{\prime}}(\dot{B}_{2,\gamma}^{\frac{5}{2\gamma^{\prime}},\frac{1}{2\gamma}})})\\
&~~~~+e^{\lambda\Psi(T)}\|(\bar{u}_{1})_{\lambda}\|_{X_{T}}\|(\bar{u}_{1})_{\lambda}\|_{\widetilde{L}_{T}^{\infty}(\dot{B}_{2,\gamma}^{\frac{1}{2\gamma^{\prime}},\frac{1}{2\gamma}})}\nonumber.
\end{align}
Combining the above estimates (\ref{hang51})-(\ref{hang1192}), one has the desired result.

Based on the above arguments, we can obtain the following a priori estimates:\\
\begin{bf}Proposition 4.4.\end{bf} Let $\lambda, k>0$ and $1<\gamma<2$. A constant C exists such that for any $T>0$ satisfying $\theta(t)\leqslant a/k$ for $t\in[0,T]$, we have
\begin{align*}&\|(\bar{u}_{1})_{\lambda}\|_{X_{T}}+\|(\Lambda_{h}^{-1}\partial_{2}\bar{u}_{2})_{\lambda}\|_{X_{T}}+\|(\Lambda_{h}^{-1}\partial_{3}\bar{u}_{3})_{\lambda}\|_{X_{T}}\\
&\leqslant \|e^{aD_{3}}u_{1}^{0}\|_{\dot{B}_{2,\gamma}^{\frac{1}{2\gamma^{\prime}},\frac{1}{2\gamma}}}
+\|\partial_{3}e^{aD_{3}}u_{3}^{0}\|_{\dot{B}_{2,\gamma}^{\frac{1}{2\gamma^{\prime}}-1,\frac{1}{2\gamma}}}+e^{\lambda\Psi(T)}
\|(\bar{u}_{1})_{\lambda}\|_{X_{T}}\|(\bar{u}_{1})_{\lambda}\|_{X_{T}}\\
&~~~~+\frac{C}{k^{\frac{1}{\gamma}}}
\|(u_{h})_{\lambda}\|_{X_{T}}+\frac{C}{\lambda^{^{\frac{1}{\gamma^{\prime}}}}}
(\|(\partial_{2}\bar{u}_{2})_{\lambda}\|_{\widetilde{L}_{T}^{\gamma}(\dot{B}_{2,\gamma}^{\frac{2}{\gamma}+\frac{1}{2\gamma^{\prime}}-1,\frac{1}{2\gamma}})}
+\|(\bar{u}_{1})_{\lambda}\|_{\widetilde{L}_{T}^{\gamma^{\prime}}(\dot{B}_{2,\gamma}^{\frac{5}{2\gamma^{\prime}},\frac{1}{2\gamma}})}).\end{align*}
Proof. Applying the Duhamel formula to the Eq.(\ref{hang111}) gives
$$\bar{u}_{3}(t)=H(t)u_{3}^{0}-\int_{0}^{t}H(t-\tau)div(u\otimes u_{3})d\tau$$
and
$$\bar{u}_{i}(t)=H(t)u_{i}^{0}-\int_{0}^{t}H(t-\tau)div_{h}(\bar{u}_{h}\bar{u}_{i})d\tau,~~i=1,2.$$
This Proposition is a direct consequence of Proposition 4.3 and Lemma 4.1.

Similar to the proof of Proposition 4.4, we get the following propositions.\\
\begin{bf}Proposition 4.5.\end{bf} Let $\lambda, k>0$ and $1<\gamma<2$. A constant C exists such that, for any $T>0$ satisfying $\theta(t)\leqslant a/k$ for $t\in[0,T]$, we have
\begin{align*}&\|(\bar{u}_{1})_{\lambda}\|_{\widetilde{L}_{T}^{\gamma}(\dot{B}_{2,\gamma}^{\frac{2}{\gamma}+\frac{1}{2\gamma^{\prime}},\frac{1}{2\gamma}})}+\|(\partial_{2}\bar{u}_{2})_{\lambda}\|_{\widetilde{L}_{T}^{\gamma}(\dot{B}_{2,\gamma}^{\frac{2}{\gamma}+\frac{1}{2\gamma^{\prime}}-1,\frac{1}{2\gamma}})}\\
&\leqslant \|e^{aD_{3}}u_{1}^{0}\|_{\dot{B}_{2,\gamma}^{\frac{1}{2\gamma^{\prime}},\frac{1}{2\gamma}}}+\|e^{aD_{3}}\partial_{2}u_{2}^{0}\|_{\dot{B}_{2,\gamma}^{\frac{1}{2\gamma^{\prime}}-1,\frac{1}{2\gamma}}}
+e^{\lambda\Psi(T)}
\|(\bar{u}_{1})_{\lambda}\|_{X_{T}}\|(\bar{u}_{1})_{\lambda}\|_{\widetilde{L}_{T}^{\gamma}(\dot{B}_{2,\gamma}^{\frac{2}{\gamma}+\frac{1}{2\gamma^{\prime}},\frac{1}{2\gamma}})}\\
&~~~~+\frac{C}{\lambda^{^{\frac{1}{\gamma^{\prime}}}}}
(\|(\partial_{2}\bar{u}_{2})_{\lambda}\|_{\widetilde{L}_{T}^{\gamma}(\dot{B}_{2,\gamma}^{\frac{2}{\gamma}+\frac{1}{2\gamma^{\prime}}-1,\frac{1}{2\gamma}})}+\|(\bar{u}_{1})_{\lambda}\|_{\widetilde{L}_{T}^{\gamma}(\dot{B}_{2,\gamma}^{\frac{2}{\gamma}+\frac{1}{2\gamma^{\prime}},\frac{1}{2\gamma}})}).\end{align*}
\begin{bf}Proposition 4.6.\end{bf} Let $\lambda, k>0$ and $1<\gamma<2$. A constant C exists such that, for any $T>0$ satisfying $\theta(t)\leqslant a/k$ for $t\in[0,T]$, we have
\begin{align*}\|(\bar{u}_{3})_{\lambda}\|_{X_{T}}
&\leqslant \|e^{aD_{3}}u_{3}^{0}\|_{\dot{B}_{2,\gamma}^{\frac{1}{2\gamma^{\prime}},\frac{1}{2\gamma}}}+(\frac{C}{\lambda^{^{\frac{1}{\gamma^{\prime}}}}}+\frac{C}{k^{\frac{1}{\gamma}}})
(\|(\bar{u}_{3})_{\lambda}\|_{X_{T}}+\|\nu_{\lambda}\|_{X_{T}})\\
&~~~~+e^{\lambda\Psi(T)}
(\|(\bar{u}_{1})_{\lambda}\|_{X_{T}}+\|\nu_{\lambda}\|_{X_{T}})(\|(\bar{u}_{3})_{\lambda}\|_{X_{T}}+\|\nu_{\lambda}\|_{X_{T}}).\end{align*}

\section{Classical analytical-type A priori estimates}

Let $X_{T}$ denote the same notation as in section 4. In this part, we need to use the fact that $e^{\Phi(t,\xi)}$ is a
sublinear function, and regularizing effect
from the analyticity.\\\\
\begin{bf}Proposition 5.1.\end{bf} Let $\lambda, k>0$ and $1<\gamma<2$. A constant C exists such that, for any $T>0$ satisfying $\theta(t)\leqslant a/k$ for $t\in[0,T]$, we have
 \begin{align}\label{hang55}
&\|u_{\Phi}\|_{L_{T}^{\gamma}(\dot{B}_{2,\gamma}^{\frac{1}{\gamma}+\frac{1}{2\gamma^{\prime}},\frac{1}{2\gamma}})\cap \widetilde{L}_{T}^{\infty}(\dot{B}_{2,\gamma}^{\frac{1}{2\gamma^{\prime}}-\frac{1}{\gamma},\frac{1}{2\gamma}})}\nonumber\\
&\lesssim
\|e^{aD_{3}}u^{0}\|_{\dot{B}_{2,\gamma}^{\frac{1}{2\gamma^{\prime}}-\frac{1}{\gamma},\frac{1}{2\gamma}}}
+\|(\bar{u}_{2})_{\Phi}\|_{L_{T}^{\gamma^{\prime}}(\dot{B}_{2,\gamma}^{\frac{5}{2\gamma^{\prime}},\frac{1}{2\gamma}})}
\|u_{\Phi}\|_{L_{T}^{\gamma}(\dot{B}_{2,\gamma}^{\frac{1}{\gamma}+\frac{1}{2\gamma^{\prime}},\frac{1}{2\gamma}})}\\
&~~~~+e^{\lambda\Psi(T)}(\|(\Lambda_{h}^{-1}\partial_{2}\bar{u}_{2})_{\lambda}\|_{X_{T}}+\|(\bar{u}_{1})_{\lambda}\|_{X_{T}}+\|\nu_{\lambda}\|_{X_{T}})\|u_{\Phi}\|_{L_{T}^{\gamma}(\dot{B}_{2,\gamma}^{\frac{1}{\gamma}+\frac{1}{2\gamma^{\prime}},\frac{1}{2\gamma}})}.
\nonumber\end{align}
\begin{bf}Proof.\end{bf} Applying the Duhamel formula to Eq.(\ref{hang113}) gives
$$u(t)=H(t)u^{0}-\int_{0}^{t}H(t-\tau)[\mathbb{P} (u\cdot\nabla u)]d\tau.$$
From Lemmas 4.1 and 4.2 and $divu=0$, it follows that
\begin{align}\label{hang557}&\|u_{\Phi}\|_{L_{T}^{\gamma}(\dot{B}_{2,\gamma}^{\frac{1}{\gamma}+\frac{1}{2\gamma^{\prime}},\frac{1}{2\gamma}})\cap \widetilde{L}_{T}^{\infty}(\dot{B}_{2,\gamma}^{\frac{1}{2\gamma^{\prime}}-\frac{1}{\gamma},\frac{1}{2\gamma}})}\nonumber\\
&\lesssim
\|e^{aD_{3}}u^{0}\|_{\dot{B}_{2,\gamma}^{\frac{1}{2\gamma^{\prime}}-\frac{1}{\gamma},\frac{1}{2\gamma}}}+\|(u_{3} div_{h}u_{h})_{\Phi}\|_{\widetilde{L}_{T}^{1}(\dot{B}_{2,\gamma}^{\frac{2}{\gamma^{\prime}}-1,\frac{1}{2\gamma}-\frac{1}{2\gamma^{\prime}}})}\\
&~~~~~+\sum_{1\leqslant i,j\leqslant3}\|(\nu_{i}u_{j})_{\Phi}\|_{\widetilde{L}_{T}^{1}(\dot{B}_{2,\gamma}^{\frac{2}{\gamma^{\prime}},\frac{1}{2\gamma}-\frac{1}{2\gamma^{\prime}}})}+\sum_{i=1,2}\|(\bar{u}_{i} u)_{\Phi}\|_{\widetilde{L}_{T}^{1}(\dot{B}_{2,\gamma}^{\frac{2}{\gamma^{\prime}},\frac{1}{2\gamma}-\frac{1}{2\gamma^{\prime}}})}.\nonumber\end{align}
Due to Lemma 3.3, we get that
\begin{align}\label{24}&\|(u_{3} div_{h}u_{h})_{\Phi}\|_{\widetilde{L}_{T}^{1}(\dot{B}_{2,\gamma}^{\frac{2}{\gamma^{\prime}}-1,\frac{1}{2\gamma}-\frac{1}{2\gamma^{\prime}}})}
+\sum_{1\leqslant i,j\leqslant3}\|(\nu_{i}u_{j})_{\Phi}\|_{\widetilde{L}_{T}^{1}(\dot{B}_{2,\gamma}^{\frac{2}{\gamma^{\prime}},\frac{1}{2\gamma}-\frac{1}{2\gamma^{\prime}}})}\nonumber\\
&\lesssim e^{\lambda\Psi(T)}(\|(\Lambda_{h}^{-1}\partial_{2}\bar{u}_{2})_{\lambda}\|_{X_{T}}+\|(\bar{u}_{1})_{\lambda}\|_{X_{T}}+\|\nu_{\lambda}\|_{X_{T}})\|u_{\Phi}\|_{L_{T}^{\gamma}(\dot{B}_{2,\gamma}^{\frac{1}{\gamma}+\frac{1}{2\gamma^{\prime}},\frac{1}{2\gamma}})}.
\end{align}
On the other hand, one can obtain that
\begin{align}\label{hang558}&\sum_{i=1,2}\|(\bar{u}_{i} u)_{\Phi}\|_{\widetilde{L}_{T}^{1}(\dot{B}_{2,\gamma}^{\frac{2}{\gamma^{\prime}},\frac{1}{2\gamma}-\frac{1}{2\gamma^{\prime}}})}\leqslant C\|(\bar{u}_{h})_{\Phi}\|_{\widetilde{L}_{T}^{\gamma^{\prime}}(\dot{B}_{2,\gamma}^{\frac{5}{2\gamma^{\prime}},\frac{1}{2\gamma}})}
\|u_{\Phi}\|_{L_{T}^{\gamma}(\dot{B}_{2,\gamma}^{\frac{1}{\gamma}+\frac{1}{2\gamma^{\prime}},\frac{1}{2\gamma}})}.
\end{align}
Plugging the above estimates into (\ref{hang557}) yields the desired results.

Similar to the proof of Proposition 4.3, one can obtain the following Proposition.\\\\
\begin{bf}Proposition 5.2.\end{bf} Let $\lambda, k>0$ and $1<\gamma<2$. A constant C exists such that, for any $T>0$ satisfying $\theta(t)\leqslant a/k$ for $t\in[0,T]$, we have
\begin{align}\label{hang551}
\|(\bar{u}_{2})_{\Phi}\|_{X_{T}}
\lesssim& \|e^{aD_{3}}u_{2}^{0}\|_{\dot{B}_{2,\gamma}^{\frac{1}{2\gamma^{\prime}},\frac{1}{2\gamma}}}+\|(\bar{u}_{1})_{\Phi}\|_{X_{T}}\|(\bar{u}_{2})_{\Phi}\|_{X_{T}}\nonumber\\
&+e^{\lambda\Psi(T)}\|(\Lambda_{h}^{-1}\partial_{2}\bar{u}_{2})_{\lambda}\|_{X_{T}}
\|(\bar{u}_{2})_{\Phi}\|_{X_{T}}\\
&+e^{\lambda\Psi(T)}\|(\partial_{2}\bar{u}_{2})_{\lambda}\|_{L_{T}^{\gamma}(\dot{B}_{2,\gamma}^{\frac{2}{\gamma}+\frac{1}{2\gamma^{\prime}}-1,\frac{1}{2\gamma}})}
\|(\bar{u}_{2})_{\Phi}\|_{\widetilde{L}_{T}^{\gamma^{\prime}}(\dot{B}_{2,\gamma}^{\frac{5}{2\gamma^{\prime}},\frac{1}{2\gamma}})}\nonumber
\end{align}
and
\begin{align}\label{hangt551}
\|(\bar{u}_{1})_{\Phi}\|_{X_{T}}
&\lesssim \|e^{aD_{3}}u_{1}^{0}\|_{\dot{B}_{2,\gamma}^{\frac{1}{2\gamma^{\prime}},\frac{1}{2\gamma}}}+\|(\bar{u}_{2})_{\Phi}\|_{L_{T}^{\gamma^{\prime}}(\dot{B}_{2,\gamma}^{\frac{5}{2\gamma^{\prime}},\frac{1}{2\gamma}})}
\|(\bar{u}_{1})_{\Phi}\|_{X_{T}}\\
&~~~~+\|(\bar{u}_{1})_{\Phi}\|_{L_{T}^{\gamma^{\prime}}(\dot{B}_{2,\gamma}^{\frac{5}{2\gamma^{\prime}},\frac{1}{2\gamma}})}\|(\bar{u}_{1})_{\Phi}\|_{X_{T}}.\nonumber
\end{align}
Proof. Applying the Duhamel formula
to the second equation
of Eq.(\ref{hang111}) gives
$$\bar{u}_{i}(t)=H(t)u_{i}^{0}-\int_{0}^{t}H(t-\tau)div_{h}(\bar{u}_{h}\bar{u}_{i})d\tau,~~i=1,2.$$
From a similar arguments as in the above (\ref{24}), one can obtain the desired results.

Next, we consider the estimates for the perturbed Navier-Stokes system (\ref{hang112}) in $\widetilde{L}_{T}^{\infty}(\dot{B}_{2,\gamma}^{\frac{1}{2\gamma^{\prime}},\frac{1}{2\gamma}})\cap \widetilde{L}_{T}^{\gamma^{\prime}}(\dot{B}_{2,\gamma}^{\frac{5}{2\gamma^{\prime}},\frac{1}{2\gamma}})$.\\\\
\begin{bf}Proposition 5.3.\end{bf} Let $\lambda, k>0$ and $1<\gamma<2$. A constant C exists such that, for any $T>0$ satisfying $\theta(t)\leqslant a/k$ for $t\in[0,T]$, we have
 \begin{align*}
\|\nu_{\lambda}\|_{X_{T}}
&\leqslant \frac{C}{k^{\frac{1}{\gamma}}}
(\|(\bar{u}_{h})_{\lambda}\|_{X_{T}}+\|\nu_{\lambda}\|_{X_{T}})+e^{\lambda\Psi(T)}
(\|(\bar{u}_{1})_{\lambda}\|_{X_{T}}+\|\nu_{\lambda}\|_{X_{T}})^{2}\\
&~~+\frac{C}{\lambda^{^{\frac{1}{\gamma^{\prime}}}}}
(\|(\partial_{2}\bar{u}_{2})_{\lambda}\|_{\widetilde{L}_{T}^{\gamma}(\dot{B}_{2,\gamma}^{\frac{2}{\gamma}+\frac{1}{2\gamma^{\prime}}-1,\frac{1}{2\gamma}})}+\|(\bar{u}_{1})_{\lambda}\|_{\widetilde{L}_{T}^{\gamma^{\prime}}(\dot{B}_{2,\gamma}^{\frac{5}{2\gamma^{\prime}},\frac{1}{2\gamma}})}
+\|\nu_{\lambda}\|_{\widetilde{L}_{T}^{\gamma^{\prime}}(\dot{B}_{2,\gamma}^{\frac{5}{2\gamma^{\prime}},\frac{1}{2\gamma}})}).
\end{align*}
\begin{bf}Proof.\end{bf} Applying the Duhamel formula to Eq.(\ref{hang112}) gives
$$\nu_{3}(t)=-\int_{0}^{t}H(t-\tau)[\mathbb{Q}(u\cdot\nabla)u]_{3}d\tau$$and
$$\nu_{i}(t)=-\int_{0}^{t}H(t-\tau)\{div_{h}(\nu_{h}\bar{u}_{i})+div_{h}(u_{h}\nu_{i})+\partial_{3}(u_{3}u_{i})+[\mathbb{Q}(u\cdot\nabla) u]_{i}\}d\tau,~~i=1,2.$$
(\ref{hang5151}) and (\ref{hang152}) implies that
\begin{align}\label{32}\|e^{-\lambda\Psi(t)+\Phi(t,D)}\mathcal{A}[\partial_{3}(u_{3}u_{i})]\|_{X_{T}}\leqslant \frac{C}{k^{\frac{1}{\gamma}}}\|(u_{h})_{\lambda}\|_{X_{T}}.\end{align}
Similar to (\ref{hang192}), one can obtain that
\begin{align}\label{ruhang192}
&\|e^{-\lambda\Psi(t)+\Phi(t,D)}\mathcal{A}[div_{h}(\nu_{h}\bar{u}_{i})+div_{h}(u_{h}\nu_{i})]\|_{X_{T}}\nonumber\\
&\lesssim\frac{C}{\lambda^{^{\frac{1}{\gamma^{\prime}}}}}
\|\nu_{\lambda}\|_{\widetilde{L}_{T}^{\gamma^{\prime}}(\dot{B}_{2,\gamma}^{\frac{5}{2\gamma^{\prime}},\frac{1}{2\gamma}})}
+e^{\lambda\Psi(T)}(\|(\bar{u}_{1})_{\lambda}\|_{X_{T}}+\|\nu_{\lambda}\|_{X_{T}})^{2}
\end{align}
Combining the above estimates (\ref{32}) and (\ref{ruhang192}) with Proposition 4.3, one can obtain the desired results.

\section{Global solution to the Navier-Stokes equation with large initial data}
The goal of this section is to present the proof of Theorem 1.
The idea of the proof is to write the solution $u(t,x)$ as follows: $$u=\bar{u}+\nu,$$ where $\bar{u}$ satisfies the transport equation (\ref{hang111}) with initial data $\bar{u}^{0}=u^{0}$. One notices that the vector field $\nu$ satisfies the perturbed Navier-Stokes system (\ref{hang112}). The proof of main result consists in studying both system. To control the some kind
of quantities on the solution, one should use the regularizing effect from the analyticity, relationship of nonlinear term between Eq.(\ref{hang111}) and Eq.(\ref{hang112}),  blow-up criterion of Eq.(\ref{hang113}) and the global well-posedness of the perturbed system (\ref{hang112}) with small data. The classical fixed point lemma, that relies on some estimates of related operators in the work space, is also used here.

\subsection{Local well-posedness of the Navier-Stokes equation (\ref{hang113}) in anisotropic space.}

In the proof of the main results, one needs to use the following local well-posedness of NS equation in anisotropic spaces. For simplicity, let $X_{T}$ denote the same notation as in section 4, and $E_{T}:=\widetilde{L}_{T}^{\infty}(\dot{B}_{2,\gamma}^{\frac{1}{2\gamma^{\prime}}-\frac{1}{\gamma},\frac{1}{2\gamma}})\cap L_{T}^{\gamma}(\dot{B}_{2,\gamma}^{\frac{1}{\gamma}+\frac{1}{2\gamma^{\prime}},\frac{1}{2\gamma}})$.\\\\
\begin{bf}Proposition 6.1.\end{bf} Let $a, k>0$, $e^{aD_{3}}u^{0}\in \dot{B}_{2,\gamma}^{\frac{1}{2\gamma^{\prime}},\frac{1}{2\gamma}}$, $1<\gamma<2$ and $\|e^{aD_{3}}u^{0}\|_{\dot{B}_{2,\gamma}^{\frac{1}{2\gamma^{\prime}}-\frac{1}{\gamma},\frac{1}{2\gamma}}}<\frac{a}{2k}$. Let $\Phi(t,\xi)$ and $\theta(t)$ be
defined by (\ref{11}) and (\ref{12}) respectively. Then there exists a $T>0$ such that the Navier-Stokes equation (\ref{hang113}) has a unique solution
$$u_{\Phi}\in X_{T}\cap E_{T}$$
and$$\theta(T)<\frac{a}{k}.$$
If $\|e^{aD_{3}}u^{0}\|_{\dot{B}_{2,\gamma}^{\frac{1}{2\gamma^{\prime}},\frac{1}{2\gamma}}}$ is sufficiently small, then the above solution is a global one, i.e., $T=\infty$.\\\\
\begin{bf}Proof.\end{bf} We consider the equivalent integral equation (\ref{21}) with respect to Eq.(\ref{hang113}). Let $T, M, N>0$ which will be chosen later. Put
\begin{eqnarray*}\mathcal {D}&=&\{u:\|u_{\Phi}\|_{X_{T}}<M;~\|u_{\Phi}\|_{E_{T}}<N;~\theta(t)<\frac{a}{k}\},\\
d(u^{\prime},u^{\prime\prime})&=&\|u_{\Phi}^{\prime}-u_{\Phi}^{\prime\prime}\|_{X_{T}\cap E_{T}}.
\end{eqnarray*}
It is easy to see that ($\mathcal {D}$,d) is a complete metric
space. Next consider the following mapping:
\begin{align*}
J:u(t)&\rightarrow e^{t\Delta}u^{0}-\int_{0}^{t}e^{(t-\tau)\Delta}\mathbb{P}(u\cdot\nabla)ud\tau\\
&:=e^{t\Delta}u^{0}-\mathcal{A}\mathbb{P}(u\cdot\nabla)u.
\end{align*}
We shall prove that there exist $T, M, N>0$ such that
$J:(\mathcal {D},d)\rightarrow(\mathcal {D},d)$ is a strict
contraction mapping.

Due to $\Phi(t,\xi)$ is a sublinear function on $\xi$, follows from Lemmas 4.2 and 3.3, one can obtain that
\begin{align}\label{hang2112}
\|\mathcal{A}[(u\cdot\nabla) u]_{\Phi}\|_{X_{T}}
&\lesssim \|u_{\Phi}\|_{\widetilde{L}_{T}^{2\gamma^{\prime}}(\dot{B}_{2,\gamma}^{\frac{3}{2\gamma^{\prime}},\frac{1}{2\gamma}})}\|u_{\Phi}\|_{X_{T}}
\end{align}
and
\begin{align}\label{hang2112222}
\|\mathcal{A}[(u\cdot\nabla) u]_{\Phi}\|_{E_{T}}
&\lesssim \|u_{\Phi}\|_{\widetilde{L}_{T}^{2\gamma^{\prime}}(\dot{B}_{2,\gamma}^{\frac{3}{2\gamma^{\prime}},\frac{1}{2\gamma}})}\|u_{\Phi}\|_{E_{T}}
\end{align}
On the other hand, in view of Lemma 4.1, one can obtain that
\begin{align*}\|e^{t\Delta}u_{\Phi}^{0}\|_{X_{T}}\lesssim\|e^{aD_{3}}u^{0}\|_{\dot{B}_{2,\gamma}^{\frac{1}{2\gamma^{\prime}},\frac{1}{2\gamma}}}\end{align*}
and
\begin{align*}\|e^{t\Delta}u_{\Phi}^{0}\|_{E_{T}}\lesssim\|e^{aD_{3}}u^{0}\|_{\dot{B}_{2,\gamma}^{\frac{1}{2\gamma^{\prime}}-\frac{1}{\gamma},\frac{1}{2\gamma}}}\end{align*}
Combining the above estimates with Proposition 5.2, one can obtain that
\begin{align*}
\|Ju_{\Phi}\|_{X_{T}}
&\leqslant \|e^{aD_{3}}u^{0}\|_{\dot{B}_{2,\gamma}^{\frac{1}{2\gamma^{\prime}},\frac{1}{2\gamma}}}+C\|u_{\Phi}\|_{\widetilde{L}_{T}^{2\gamma^{\prime}}(\dot{B}_{2,\gamma}^{\frac{3}{2\gamma^{\prime}},\frac{1}{2\gamma}})}
\|u_{\Phi}\|_{X_{T}};
\end{align*}
\begin{align*}
\|Ju_{\Phi}\|_{E_{T}}
&\leqslant \|e^{aD_{3}}u^{0}\|_{\dot{B}_{2,\gamma}^{\frac{1}{2\gamma^{\prime}}-\frac{1}{\gamma},\frac{1}{2\gamma}}}+C\|u_{\Phi}\|_{\widetilde{L}_{T}^{2\gamma^{\prime}}(\dot{B}_{2,\gamma}^{\frac{3}{2\gamma^{\prime}},\frac{1}{2\gamma}})}
\|u_{\Phi}\|_{E_{T}}
\end{align*}
and\begin{align*}
\|Ju_{\Phi}^{\prime}-Ju_{\Phi}^{\prime\prime}\|_{X_{T}\cap E_{T}}
&\lesssim
(\|u_{\Phi}^{\prime}\|_{\widetilde{L}_{T}^{2\gamma^{\prime}}(\dot{B}_{2,\gamma}^{\frac{3}{2\gamma^{\prime}},\frac{1}{2\gamma}})}
+\|u_{\Phi}^{\prime\prime}\|_{\widetilde{L}_{T}^{2\gamma^{\prime}}(\dot{B}_{2,\gamma}^{\frac{3}{2\gamma^{\prime}},\frac{1}{2\gamma}})})\|u_{\Phi}^{\prime}-u_{\Phi}^{\prime\prime}\|_{X_{T}\cap E_{T}}.
\end{align*}
Put $M=2\|e^{aD_{3}}u^{0}\|_{\dot{B}_{2,\gamma}^{\frac{1}{2\gamma^{\prime}},\frac{1}{2\gamma}}}$ and $N=2\|e^{aD_{3}}u^{0}\|_{\dot{B}_{2,\gamma}^{\frac{1}{2\gamma^{\prime}}-\frac{1}{\gamma},\frac{1}{2\gamma}}}$. Taking $T>0$ satisfying $$C\|u_{\Phi}\|_{\widetilde{L}_{T}^{2\gamma^{\prime}}(\dot{B}_{2,\gamma}^{\frac{3}{2\gamma^{\prime}},\frac{1}{2\gamma}})}<\frac{1}{4}.$$
Then one has that
$\mathcal{J}:(\mathcal{D},d)\rightarrow(\mathcal{D},d)$ is a contraction mapping. So, there exists a $u\in\mathcal{D}$ satisfying NS equation (\ref{hang113}).
For the global wellposedness with small data, it suffices to take
$T=\infty$ in ($\mathcal {D}$,d).
Furthermore, relying on a continuation argument, one can also obtain the local well-posedness of system (\ref{hang111}) similarly.

\subsection{Proof of Theorem 1.}

In this subsection, we will use Propositions in section 4 and section 5 to complete the proof of Theorem 1.\\\\
\begin{bf}End of the proof of Theorem 1.\end{bf}
Now we are ready to prove Theorem 1. The idea, as presented in the introduction, is to write
$$u=\bar{u}+\nu,$$ where $\nu$ satisfies the perturbed system (\ref{hang112}) with $\nu^{0}=0$.

We use a continuation argument. For any
$\eta$, let us define
\begin{align*}T_{\lambda, k}:=\left\{T:
  \begin{array}{lll}
  \theta(T)^{\frac{1}{\gamma}}\leqslant \frac{3}{2}C\exp\{4C\|e^{aD_{3}}u_{2}^{0}\|_{\dot{B}_{2,\gamma}^{\frac{1}{2\gamma^{\prime}},\frac{1}{2\gamma}}}^{\gamma^{\prime}}\}\|e^{aD_{3}}u^{0}\|_{\dot{B}_{2,\gamma}^{\frac{1}{2\gamma^{\prime}}-\frac{1}{\gamma},\frac{1}{2\gamma}}};\\
  \|(\bar{u}_{3})_{\lambda}\|_{L_{T}^{\gamma^{\prime}}(\dot{B}_{2,\gamma}^{\frac{5}{2\gamma^{\prime}},\frac{1}{2\gamma}})}
  \leqslant 4C\|e^{aD_{3}}u_{3}^{0}\|_{\dot{B}_{2,\gamma}^{\frac{1}{2\gamma^{\prime}},\frac{1}{2\gamma}}}+4C\eta;\\
  \|(\bar{u}_{1})_{\Phi}\|_{L_{T}^{\gamma^{\prime}}(\dot{B}_{2,\gamma}^{\frac{5}{2\gamma^{\prime}},\frac{1}{2\gamma}})}\leqslant4C\eta;~\Psi(T)^{\frac{1}{\gamma^{\prime}}}\leqslant4C\|e^{aD_{3}}u_{2}^{0}\|_{\dot{B}_{2,\gamma}^{\frac{1}{2\gamma^{\prime}},\frac{1}{2\gamma}}};\\
  \exp\{4C\lambda\|e^{aD_{3}}u_{2}^{0}\|_{\dot{B}_{2,\gamma}^{\frac{1}{2\gamma^{\prime}},\frac{1}{2\gamma}}}^{\gamma^{\prime}}\}(\|(\bar{u}_{1})_{\lambda}\|_{X_{T}}+\|(\Lambda_{h}^{-1}\partial_{2}\bar{u}_{2})_{\lambda}\|_{X_{T}}\\
  +\|(\Lambda_{h}^{-1}\partial_{3}\bar{u}_{3})_{\lambda}\|_{X_{T}}+\|\nu_{\lambda}\|_{X_{T}}+\|(\partial_{2}\bar{u}_{2})_{\lambda}\|_{\widetilde{L}_{T}^{\gamma}(\dot{B}_{2,\gamma}^{\frac{2}{\gamma}+\frac{1}{2\gamma^{\prime}}-1,\frac{1}{2\gamma}})})\leqslant4C\eta
  \end{array}
  \right\}.
\end{align*}
Similar to the argument in \cite{Chemin4} (the local well-posedness of Eq.(\ref{hang113}) will be used here), $T_{\lambda, k}$ is of the form $[0,T^{*})$ for some positive $T^{*}$. Hence,
it suffices to prove that $T^{*}=\infty$. In order to use Propositions in sections 4 and 5, we need to
assume that \begin{align*}\theta(t)\leqslant\frac{a}{k},\end{align*} which leads to the condition
\begin{align}\label{30}\frac{3}{2}C\exp\{4C\|e^{aD_{3}}u_{2}^{0}\|_{\dot{B}_{2,\gamma}^{\frac{1}{2\gamma^{\prime}},\frac{1}{2\gamma}}}^{\gamma^{\prime}}\}\|e^{aD_{3}}u^{0}\|_{\dot{B}_{2,\gamma}^{\frac{1}{2\gamma^{\prime}}-\frac{1}{\gamma},\frac{1}{2\gamma}}}\leqslant (\frac{a}{k})^{\frac{1}{\gamma}}.\end{align}
 From Proposition 5.2, it follows that, for all $T\in T_{\lambda, k}$,
 \begin{align}\label{25}(i),~~&\|(\bar{u}_{2})_{\Phi}\|_{L_{T}^{\gamma^{\prime}}(\dot{B}_{2,\gamma}^{\frac{5}{2\gamma^{\prime}},\frac{1}{2\gamma}})}\leqslant 2C\|e^{aD_{3}}u_{2}^{0}\|_{\dot{B}_{2,\gamma}^{\frac{1}{2\gamma^{\prime}},\frac{1}{2\gamma}}}.
\end{align}
Splitting [0, T] into m subintervals like as [0, $T_{1}$], [$T_{1}$, $T_{2}$] and so on, such that
$$C\|(\bar{u}_{2})_{\Phi}\|_{L^{\gamma^{\prime}}([T_{i},T_{i+1});\dot{B}_{2,\gamma}^{\frac{5}{2\gamma^{\prime}},\frac{1}{2\gamma}})}<\frac{1}{4},~~for~i=0,2,......m-1,$$
with $T_{0}=0$. Arguing similarly as in deriving (\ref{hang55}) and (\ref{hangt551}), one can get for all $t\in[T_{i}, T_{i+1})$,
\begin{align*}&\|u_{\Phi}\|_{L^{\gamma}([T_{i}, T_{i+1});\dot{B}_{2,\gamma}^{\frac{1}{\gamma}+\frac{1}{2\gamma^{\prime}},\frac{1}{2\gamma}})\cap \widetilde{L}^{\infty}([T_{i}, T_{i+1});\dot{B}_{2,\gamma}^{\frac{1}{2\gamma^{\prime}}-\frac{1}{\gamma},\frac{1}{2\gamma}})}\\
 &\leqslant C\|e^{\Phi(T_{i},D)}u(T_{i})\|_{\dot{B}_{2,\gamma}^{\frac{1}{2\gamma^{\prime}}-\frac{1}{\gamma},\frac{1}{2\gamma}}}+4C\eta\|u_{\Phi}\|_{L^{\gamma}([T_{i}, T_{i+1});\dot{B}_{2,\gamma}^{\frac{1}{\gamma}+\frac{1}{2\gamma^{\prime}},\frac{1}{2\gamma}})}\\
 &~~~~+C\|(\bar{u}_{2})_{\Phi}\|_{L^{\gamma^{\prime}}([T_{i}, T_{i+1});\dot{B}_{2,\gamma}^{\frac{5}{2\gamma^{\prime}},\frac{1}{2\gamma}})}\|u_{\Phi}\|_{L^{\gamma}([T_{i}, T_{i+1});\dot{B}_{2,\gamma}^{\frac{1}{\gamma}+\frac{1}{2\gamma^{\prime}},\frac{1}{2\gamma}})},
\end{align*}
and
\begin{align*}&\|(\bar{u}_{1})_{\Phi}\|_{L^{\gamma^{\prime}}([T_{i}, T_{i+1});\dot{B}_{2,\gamma}^{\frac{5}{2\gamma^{\prime}},\frac{1}{2\gamma}})\cap \widetilde{L}^{\infty}([T_{i}, T_{i+1});\dot{B}_{2,\gamma}^{\frac{1}{2\gamma^{\prime}},\frac{1}{2\gamma}})}\\
&\leqslant C\|e^{\Phi(T_{i},D)}u_{1}(T_{i})\|_{\dot{B}_{2,\gamma}^{\frac{1}{2\gamma^{\prime}},\frac{1}{2\gamma}}}+4C\eta \|(\bar{u}_{1})_{\Phi}\|_{L^{\gamma^{\prime}}([T_{i}, T_{i+1});\dot{B}_{2,\gamma}^{\frac{5}{2\gamma^{\prime}},\frac{1}{2\gamma}})\cap \widetilde{L}^{\infty}([T_{i}, T_{i+1});\dot{B}_{2,\gamma}^{\frac{1}{2\gamma^{\prime}},\frac{1}{2\gamma}})}\\
&~~~~+C\|(\bar{u}_{2})_{\Phi}\|_{L^{\gamma^{\prime}}([T_{i}, T_{i+1});\dot{B}_{2,\gamma}^{\frac{5}{2\gamma^{\prime}},\frac{1}{2\gamma}})} \|(\bar{u}_{1})_{\Phi}\|_{L^{\gamma^{\prime}}([T_{i}, T_{i+1});\dot{B}_{2,\gamma}^{\frac{5}{2\gamma^{\prime}},\frac{1}{2\gamma}})\cap \widetilde{L}^{\infty}([T_{i}, T_{i+1});\dot{B}_{2,\gamma}^{\frac{1}{2\gamma^{\prime}},\frac{1}{2\gamma}})}.
\end{align*}
Since the number of such subintervals is $m\thickapprox C\|e^{aD_{3}}u_{2}^{0}\|_{\dot{B}_{2,\gamma}^{\frac{1}{2\gamma^{\prime}},\frac{1}{2\gamma}}}^{\gamma^{\prime}}$, by a standard
induction argument, one can readily conclude that for all $T\in T_{\lambda, k}$,
 \begin{align*}(ii),~~&\theta(T)^{\frac{1}{\gamma}}\leqslant C\exp\{4C\|e^{aD_{3}}u_{2}^{0}\|_{\dot{B}_{2,\gamma}^{\frac{1}{2\gamma^{\prime}},\frac{1}{2\gamma}}}^{\gamma^{\prime}}\}\|e^{aD_{3}}u^{0}\|_{\dot{B}_{2,\gamma}^{\frac{1}{2\gamma^{\prime}}-\frac{1}{\gamma},\frac{1}{2\gamma}}},\end{align*}and
  \begin{align*}
(iii),~~&\|(\bar{u}_{1})_{\Phi}\|_{L_{T}^{\gamma^{\prime}}(\dot{B}_{2,\gamma}^{\frac{5}{2\gamma^{\prime}},\frac{1}{2\gamma}})}\leqslant C\exp\{4C\|e^{aD_{3}}u_{2}^{0}\|_{\dot{B}_{2,\gamma}^{\frac{1}{2\gamma^{\prime}},\frac{1}{2\gamma}}}^{\gamma^{\prime}}\}\|e^{aD_{3}}u_{1}^{0}\|_{\dot{B}_{2,\gamma}^{\frac{1}{2\gamma^{\prime}},\frac{1}{2\gamma}}}.\end{align*}
Propositions 4.5 and 4.6 implies that for all $T\in T_{\lambda, k}$,
\begin{align*}&(iv),~~\|(\partial_{2}\bar{u}_{2})_{\lambda}\|_{\widetilde{L}_{T}^{\gamma}(\dot{B}_{2,\gamma}^{\frac{2}{\gamma}+\frac{1}{2\gamma^{\prime}}-1,\frac{1}{2\gamma}})}\leqslant 2C(\|e^{aD_{3}}u_{1}^{0}\|_{\dot{B}_{2,\gamma}^{\frac{1}{2\gamma^{\prime}},\frac{1}{2\gamma}}}+\|e^{aD_{3}}\partial_{2}u_{2}^{0}\|_{\dot{B}_{2,\gamma}^{\frac{1}{2\gamma^{\prime}}-1,\frac{1}{2\gamma}}})
\end{align*}and
\begin{align*}(v),~~\|(\bar{u}_{3})_{\lambda}\|_{L_{T}^{\gamma^{\prime}}(\dot{B}_{2,\gamma}^{\frac{5}{2\gamma^{\prime}},\frac{1}{2\gamma}})}
&\leqslant 2C\|e^{aD_{3}}u_{3}^{0}\|_{\dot{B}_{2,\gamma}^{\frac{1}{2\gamma^{\prime}},\frac{1}{2\gamma}}}+2C\eta.\end{align*}
Thanks to Propositions 4.4 and 5.3, one can obtain that for all $T\in T_{\lambda, k}$,
\begin{align*}(vi),~~&\exp\{4C\lambda\|e^{aD_{3}}u_{2}^{0}\|_{\dot{B}_{2,\gamma}^{\frac{1}{2\gamma^{\prime}},\frac{1}{2\gamma}}}^{\gamma^{\prime}}\}\\
&~~~~\times(\|(\bar{u}_{1})_{\lambda}\|_{X_{T}}
+\|(\Lambda_{h}^{-1}\partial_{2}\bar{u}_{2})_{\lambda}\|_{X_{T}}+\|(\Lambda_{h}^{-1}\partial_{3}\bar{u}_{3})_{\lambda}\|_{X_{T}})\\
&\leqslant C\exp\{4C\lambda\|e^{aD_{3}}u_{2}^{0}\|_{\dot{B}_{2,\gamma}^{\frac{1}{2\gamma^{\prime}},\frac{1}{2\gamma}}}^{\gamma^{\prime}}\}(\|e^{aD_{3}}u_{1}^{0}\|_{\dot{B}_{2,\gamma}^{\frac{1}{2\gamma^{\prime}},\frac{1}{2\gamma}}}
+\|\partial_{3}e^{aD_{3}}u_{3}^{0}\|_{\dot{B}_{2,\gamma}^{\frac{1}{2\gamma^{\prime}}-1,\frac{1}{2\gamma}}})\\
&~~~~+\frac{C}{k^{\frac{1}{\gamma}}}\exp\{4C\lambda\|e^{aD_{3}}u_{2}^{0}\|_{\dot{B}_{2,\gamma}^{\frac{1}{2\gamma^{\prime}},\frac{1}{2\gamma}}}^{\gamma^{\prime}}\}
\|(\bar{u}_{2})_{\lambda}\|_{X_{T}}+4C(\frac{1}{\lambda^{^{\frac{1}{\gamma^{\prime}}}}}+\frac{1}{k^{\frac{1}{\gamma}}})\eta+C\eta^{2},\\
(vii),~~&\exp\{4C\lambda\|e^{aD_{3}}u_{2}^{0}\|_{\dot{B}_{2,\gamma}^{\frac{1}{2\gamma^{\prime}},\frac{1}{2\gamma}}}^{\gamma^{\prime}}\}\|\nu_{\lambda}\|_{X_{T}}\\
&\leqslant\frac{C}{k^{\frac{1}{\gamma}}}
\exp\{4C\lambda\|e^{aD_{3}}u_{2}^{0}\|_{\dot{B}_{2,\gamma}^{\frac{1}{2\gamma^{\prime}},\frac{1}{2\gamma}}}^{\gamma^{\prime}}\}\|(\bar{u}_{2})_{\lambda}\|_{X_{T}}
+4C(\frac{1}{\lambda^{^{\frac{1}{\gamma^{\prime}}}}}+\frac{1}{k^{\frac{1}{\gamma}}})\eta+C\eta^{2}.\end{align*}
 We first choose $\lambda$ large enough
that $$\frac{1}{\lambda^{\frac{1}{\gamma^{\prime}}}}\ll1.$$
Moreover, let $a\geqslant1$, in addition to (\ref{30}), one can further choose $k$ large enough such that
 $$(\frac{a}{k})^{\frac{1}{\gamma}}\leqslant 2C\exp\{4C\|e^{aD_{3}}u_{2}^{0}\|_{\dot{B}_{2,\gamma}^{\frac{1}{2\gamma^{\prime}},\frac{1}{2\gamma}}}^{\gamma^{\prime}}\}\|e^{aD_{3}}u^{0}\|_{\dot{B}_{2,\gamma}^{\frac{1}{2\gamma^{\prime}}-\frac{1}{\gamma},\frac{1}{2\gamma}}}.$$
With this choice of $k$, we infer from the condition (\ref{hang3332}) in Theorem 1 that
\begin{align}\label{33}\frac{C}{k^{\frac{1}{\gamma}}}
\exp\{4C\lambda\|e^{aD_{3}}u_{2}^{0}\|_{\dot{B}_{2,\gamma}^{\frac{1}{2\gamma^{\prime}},\frac{1}{2\gamma}}}^{\gamma^{\prime}}\}\|e^{aD_{3}}u_{2}^{0}\|_{\dot{B}_{2,\gamma}^{\frac{1}{2\gamma^{\prime}},\frac{1}{2\gamma}}}<\eta.\end{align}
Let $\eta\ll1$,
 combining (\ref{33}) with the above estimates (i)-(vii), we obtain that
\begin{align*}
&\theta(T)^{\frac{1}{\gamma}}< \frac{3}{2}C\exp\{4C\|e^{aD_{3}}u_{2}^{0}\|_{\dot{B}_{2,\gamma}^{\frac{1}{2\gamma^{\prime}},\frac{1}{2\gamma}}}^{\gamma^{\prime}}\}\|e^{aD_{3}}u^{0}\|_{\dot{B}_{2,\gamma}^{\frac{1}{2\gamma^{\prime}}-\frac{1}{\gamma},\frac{1}{2\gamma}}};\\
&\|(\bar{u}_{3})_{\lambda}\|_{L_{T}^{\gamma^{\prime}}(\dot{B}_{2,\gamma}^{\frac{5}{2\gamma^{\prime}},\frac{1}{2\gamma}})}
< 4C\|e^{aD_{3}}u_{3}^{0}\|_{\dot{B}_{2,\gamma}^{\frac{1}{2\gamma^{\prime}},\frac{1}{2\gamma}}}+4C\eta;\\
&\|(\bar{u}_{1})_{\Phi}\|_{L_{T}^{\gamma^{\prime}}(\dot{B}_{2,\gamma}^{\frac{5}{2\gamma^{\prime}},\frac{1}{2\gamma}})}<4C\eta;~~\Psi(T)^{\frac{1}{\gamma^{\prime}}}<4C\|e^{aD_{3}}u_{2}^{0}\|_{\dot{B}_{2,\gamma}^{\frac{1}{2\gamma^{\prime}},\frac{1}{2\gamma}}};\\
  &\exp\{4C\lambda\|e^{aD_{3}}u_{2}^{0}\|_{\dot{B}_{2,\gamma}^{\frac{1}{2\gamma^{\prime}},\frac{1}{2\gamma}}}^{\gamma^{\prime}}\}(\|(\partial_{2}\bar{u}_{2})_{\lambda}\|_{\widetilde{L}_{T}^{\gamma}(\dot{B}_{2,\gamma}^{\frac{2}{\gamma}+\frac{1}{2\gamma^{\prime}}-1,\frac{1}{2\gamma}})}+\|(\bar{u}_{1})_{\lambda}\|_{X_{T}}\\
  &~~~~~~~~~~~~~+\|(\Lambda_{h}^{-1}\partial_{2}\bar{u}_{2})_{\lambda}\|_{X_{T}}+\|(\Lambda_{h}^{-1}\partial_{3}\bar{u}_{3})_{\lambda}\|_{X_{T}}+\|\nu_{\lambda}\|_{X_{T}})<4C\eta,
\end{align*}
which ensures that $T^{*}=\infty$, and thus we conclude the proof of Theorem 1.

\subsection{Proof of Theorem 2.}

In this section, we shall prove that the scaling of the initial data claimed in Theorem 2, satisfies the assumptions in the Theorem 1, namely the assumption (\ref{hang3332}), in which case the global
wellposedness will follow as a consequence of that Theorem. In the proof, one needs to use the anisotropic
spaces and their propositions.

First, similar to the arguments as in \cite{Chemin3}, one can obtain the following Proposition.\\\\
\begin{bf}Proposition 6.2.\end{bf} Let $(f,g,h)$ be in $\mathcal{S}(\mathbb{R}^{3})$, $0\leqslant\alpha\leqslant\beta<1$. Let us define $$\psi^{\varepsilon}(x):=f(\varepsilon^{\alpha}x_{1})g(\varepsilon^{\beta} x_{2})h(\varepsilon x_{3}).$$ We have,
if $\varepsilon$ is small enough, \begin{align}\label{13}\|\psi^{\varepsilon}\|_{\dot{B}_{\infty,\infty}^{-1}(\mathbb{R}^{3})}&\geqslant\frac{1}{4}\varepsilon^{-\alpha}\|f\|_{\dot{B}_{\infty,\infty}^{-1}(\mathbb{R})}\|g\|_{L^{\infty}(\mathbb{R})}\|h\|_{L^{\infty}(\mathbb{R})},~~for~~\alpha<\beta,\\
\|\psi^{\varepsilon}\|_{\dot{B}_{\infty,\infty}^{-1}(\mathbb{R}^{3})}&\geqslant\frac{1}{4}\varepsilon^{-\alpha}\|fg\|_{\dot{B}_{\infty,\infty}^{-1}(\mathbb{R}^{2})}\|h\|_{L^{\infty}(\mathbb{R})},~~for~~\alpha=\beta.~~~~~~~~~\end{align}
\begin{bf}Proof.\end{bf} By the definition of $\|\cdot\|_{\dot{B}_{\infty,\infty}^{-1}(\mathbb{R}^{3})}$ given by (d) in Proposition 2.2, we have to bound from below the
quantity $\|e^{t\Delta}\psi^{\varepsilon}\|_{L^{\infty}(\mathbb{R}^{3})}$. Write
$$e^{t\Delta}\psi^{\varepsilon}(t,x)=(e^{t\partial_{1}^{2}}f)(\varepsilon^{2\alpha}t,\varepsilon^{\alpha} x_{1})(e^{t\partial_{2}^{2}}g)(\varepsilon^{2\beta}t,\varepsilon^{\beta} x_{2})(e^{t\partial_{3}^{2}}h)(\varepsilon^{2}t,\varepsilon x_{3}).$$
Then we have
\begin{align*}&\sup\limits_{t>0}t^{\frac{1}{2}}\|e^{t\Delta}\psi^{\varepsilon}\|_{L^{\infty}(\mathbb{R}^{3})}\\
&=\varepsilon^{-\alpha}\sup\limits_{t>0}t^{\frac{1}{2}}\|(e^{t\partial_{1}^{2}}f)(t,\varepsilon^{\alpha} x_{1})(e^{t\partial_{2}^{2}}g)(\varepsilon^{2(\beta-\alpha)}t,\varepsilon^{\beta} x_{2})(e^{t\partial_{3}^{2}}h)(\varepsilon^{2-2\alpha}t,\varepsilon x_{3})\|_{L^{\infty}(\mathbb{R}^{3})}.\end{align*}
On the other hand, let us consider a positive time $t_{0}$ such that
$$t_{0}^{\frac{1}{2}}\|e^{t_{0}\partial_{1}^{2}}f\|_{L^{\infty}(\mathbb{R})}\geqslant\frac{1}{2}\|f\|_{\dot{B}_{\infty,\infty}^{-1}(\mathbb{R})}.$$ Then, we have
\begin{align*}
&t_{0}^{\frac{1}{2}}\|e^{t_{0}\partial_{1}^{2}}f\|_{L^{\infty}(\mathbb{R})}\|(e^{t\partial_{2}^{2}}g)(\varepsilon^{2(\beta-\alpha)}t_{0},\varepsilon^{\beta} x_{2})\|_{L^{\infty}(\mathbb{R})}\|(e^{t\partial_{3}^{2}}h)(\varepsilon^{(2-2\alpha)}t_{0},\varepsilon x_{3})\|_{L^{\infty}(\mathbb{R})}\\
&\geqslant\frac{1}{2}\|f\|_{\dot{B}_{\infty,\infty}^{-1}(\mathbb{R})}\|(e^{t\partial_{2}^{2}}g)(\varepsilon^{2(\beta-\alpha)}t_{0},\varepsilon^{\beta} x_{2})\|_{L^{\infty}(\mathbb{R})}\|(e^{t\partial_{3}^{2}}h)(\varepsilon^{(2-2\alpha)}t_{0},\varepsilon x_{3})\|_{L^{\infty}(\mathbb{R})}.
\end{align*}
As, for $\alpha<\beta$, $\lim\limits_{\varepsilon\rightarrow0}e^{\varepsilon^{2(\beta-\alpha)}t_{0}\partial_{2}^{2}}g=g$ and $\lim\limits_{\varepsilon\rightarrow0}e^{\varepsilon^{(2-2\alpha)}t_{0}\partial_{3}^{2}}h=h$ in $L^{\infty}(\mathbb{R})$, (\ref{13}) is proved. For the case $\alpha=\beta$, one can prove it similarly.\\\\
\begin{bf}Proof of Theorem 2.\end{bf}
As stated in the above Proposition 6.2, one can see that $\|u_{\varepsilon,\alpha,\beta,\sigma}^{0}\|_{\dot{B}_{\infty,\infty}^{-1}}\thicksim\varepsilon^{-(\alpha\wedge\beta)-\sigma}$. Next, by using the results in Theorem 1, we will prove that $u_{\varepsilon,\alpha,\beta,\sigma}^{0}$ generates a unique global solution to Eq.(\ref{hang113}).

First, we consider the following initial data:
 \begin{align*}\widetilde{u}_{\varepsilon,\alpha,\beta,\sigma}^{0}(x)&=(\varepsilon^{-(\alpha+\sigma)} \omega_{1}^{0}(\varepsilon^{\alpha-1}x_{1},\varepsilon^{\beta-1} x_{2},x_{3}),\varepsilon^{-(\beta+\sigma)} \omega_{2}^{0}(\varepsilon^{\alpha-1}x_{1},\varepsilon^{\beta-1} x_{2},x_{3}),\\
&~~~~~~~\varepsilon^{-(1+\sigma)}\omega_{3}^{0}(\varepsilon^{\alpha-1}x_{1},\varepsilon^{\beta-1} x_{2},x_{3})),\end{align*}
with $1-\frac{3}{2}\alpha-\frac{1}{2}\beta>\sigma$ and $1-\frac{1}{2}\alpha-\frac{3}{2}\beta>\sigma$.
It is easy to see that there exists a constant C, which independent of $\varepsilon$, such that
$$\|e^{aD_{3}}(\widetilde{u}_{\varepsilon,\alpha,\beta,\sigma}^{0})_{1}\|_{\dot{B}_{2,\gamma}^{\frac{1}{2\gamma^{\prime}},\frac{1}{2\gamma}}}\leqslant C\varepsilon^{1-\frac{3}{2}\alpha-\frac{1}{2}\beta-\sigma+\frac{1}{2\gamma^{\prime}}(\alpha\wedge\beta-1)}\|e^{aD_{3}}\omega_{1}^{0}\|_{\dot{B}_{2,\gamma}^{\frac{1}{2\gamma^{\prime}},\frac{1}{2\gamma}}}$$and
$$\|e^{aD_{3}}(\widetilde{u}_{\varepsilon,\alpha,\beta,\sigma}^{0})_{2}\|_{\dot{B}_{2,\gamma}^{\frac{1}{2\gamma^{\prime}},\frac{1}{2\gamma}}}\leqslant C\varepsilon^{1-\frac{1}{2}\alpha-\frac{3}{2}\beta-\sigma+\frac{1}{2\gamma^{\prime}}(\alpha\wedge\beta-1)}\|e^{aD_{3}}\omega_{2}^{0}\|_{\dot{B}_{2,\gamma}^{\frac{1}{2\gamma^{\prime}},\frac{1}{2\gamma}}}.$$
Similarly, one can obtain the estimates for the initial data in $\dot{B}_{2,\gamma}^{\frac{1}{2\gamma^{\prime}}-\frac{1}{\gamma},\frac{1}{2\gamma}}$ that
$$\|e^{aD_{3}}(\widetilde{u}_{\varepsilon,\alpha,\beta,\sigma}^{0})_{h}\|_{\dot{B}_{2,\gamma}^{\frac{1}{2\gamma^{\prime}}-\frac{1}{\gamma},\frac{1}{2\gamma}}}\leqslant C\varepsilon^{\frac{(2-\alpha-\beta)(3+\gamma)}{4\gamma}-(\alpha\vee\beta+\sigma)}\|e^{aD_{3}}\omega_{h}^{0}\|_{\dot{B}_{\frac{4\gamma}{3+\gamma},\gamma}^{0,\frac{5-\gamma}{4\gamma}}}$$
and
$$\|e^{aD_{3}}(\widetilde{u}_{\varepsilon,\alpha,\beta,\sigma}^{0})_{3}\|_{\dot{B}_{2,\gamma}^{\frac{1}{2\gamma^{\prime}}-\frac{1}{\gamma},\frac{1}{2\gamma}}}\leqslant C\varepsilon^{[\frac{(2-\alpha-\beta)(3+\gamma)}{4\gamma}-(1+\sigma)]}\|e^{aD_{3}}\omega_{3}^{0}\|_{\dot{B}_{\frac{4\gamma}{3+\gamma},\gamma}^{0,\frac{5-\gamma}{4\gamma}}}.$$
Furthermore, one can infer from embedding proposition for anisotropic homogeneous Besov spaces that
$$\|e^{aD_{3}}\partial_{3}(\widetilde{u}_{\varepsilon,\alpha,\beta,\sigma}^{0})_{2}\|_{\dot{B}_{2,\gamma}^{\frac{1}{2\gamma^{\prime}}-1,\frac{1}{2\gamma}}}\leqslant C\varepsilon^{[\frac{(2-\alpha-\beta)(1+3\gamma)}{4\gamma}-(\beta+\sigma)]}\|e^{aD_{3}}\partial_{3}\omega_{2}^{0}(x)\|_{\dot{B}_{\frac{4\gamma}{1+3\gamma},\gamma}^{0,\frac{3+\gamma}{4\gamma}}}$$
and
$$\|e^{aD_{3}}\partial_{3}(\widetilde{u}_{\varepsilon,\alpha,\beta,\sigma}^{0})_{3}\|_{\dot{B}_{2,\gamma}^{\frac{1}{2\gamma^{\prime}}-1,\frac{1}{2\gamma}}}\leqslant C\varepsilon^{[\frac{(2-\alpha-\beta)(1+3\gamma)}{4\gamma}-(1+\sigma)]}\|e^{aD_{3}}\partial_{3}\omega_{3}^{0}(x)\|_{\dot{B}_{\frac{4\gamma}{1+3\gamma},\gamma}^{0,\frac{3+\gamma}{4\gamma}}}.$$
Based on the above observation, selecting $\gamma$ close to 1, such that $$1-\frac{3}{2}\alpha-\frac{1}{2}\beta-\sigma+\frac{1}{2\gamma^{\prime}}(\alpha\wedge\beta-1)>0,$$
$$1-\frac{1}{2}\alpha-\frac{3}{2}\beta-\sigma+\frac{1}{2\gamma^{\prime}}(\alpha\wedge\beta-1)>0,$$
$$\frac{(2-\alpha-\beta)(3+\gamma)}{4\gamma}-(1+\sigma)>0$$ and $$\frac{(2-\alpha-\beta)(1+3\gamma)}{4\gamma}-(1+\sigma)>0,$$ one can obtain that the initial data $\widetilde{u}_{\varepsilon,\alpha,\beta,\sigma}^{0}(x)$ satisfies the assumptions (\ref{hang3332}) in Theorem 1. So, $\widetilde{u}_{\varepsilon,\alpha,\beta,\sigma}^{0}(x)$ generates a unique global solution to Eq.(\ref{hang113}).
From scaling invariant (\ref{18}), and taking $\lambda=\varepsilon$, one can obtain the desired results.

\section{Appendix.}
In this section, we will list several useful Propositions.\\\\
\begin{bf}Proposition A.\end{bf}
Let $1\leqslant r\leqslant p\leqslant\infty$, and
$k=0,1$. Then
$$\|\nabla^{k}H(t)f\|_{L^{p}(\mathbb{R}^{d})}\lesssim
t^{-\frac{k}{2}-\frac{1}{2}(\frac{d}{r}-\frac{d}{p})}\|f\|_{L^{r}(\mathbb{R}^{d})},~t>0.$$
Proof. One can refer to \cite{Wang} for the desired result.

 As a generalization, from a similar argument as in Proposition A, one can show that for any $s\geqslant0$, $1\leqslant q\leqslant\infty$ and $1\leqslant r\leqslant p\leqslant\infty$,
 $$\|\mathcal{F}_{x_{3}}(-\Delta_{h})^{\frac{s}{2}}H(t)f\|_{L_{x_{3}}^{q}L_{x_{h}}^{p}(\mathbb{R}^{3})}\lesssim
t^{-\frac{s}{2}-\frac{1}{2}(\frac{2}{r}-\frac{2}{p})}\|\mathcal{F}_{x_{3}}f\|_{L_{x_{3}}^{q}L_{x_{h}}^{r}(\mathbb{R}^{3})};$$
$$\|(-\Delta)^{\frac{s}{2}}H(t)f\|_{L^{p}(\mathbb{R}^{d})}\lesssim
t^{-\frac{s}{2}-\frac{1}{2}(\frac{d}{r}-\frac{d}{p})}\|f\|_{L^{r}(\mathbb{R}^{d})},$$
where $(-\Delta)^{\frac{s}{2}}f=\mathcal{F}^{-1}|\xi|^{s}\mathcal{F}f$.

Up to now, the singular integration is still a difficult problem in the theory of
harmonic analysis. Hardy-Littlewood-Sobolev inequality is a fundamental
tool in this subject. Let $0<\alpha<d$, $$I_{\alpha}f=\int_{\mathbb{R}^{n}}\frac{f(y)}{|x-y|^{d-\alpha}}dy.$$
\begin{bf}Proposition B.\end{bf} \cite{Stein} Let $1<p,q<\infty$ and $0<\alpha<d$ satisfy $\frac{1}{p}=\frac{1}{q}+\frac{\alpha}{d}$. Then we have
$$\|I_{\alpha}f\|_{L^{q}}\leqslant \|f\|_{L^{p}}.$$\\
\begin{bf}Definition C.\end{bf} \cite{Chemin1} Let A and B be two positive real numbers and let p in $(3,+\infty]$. We define
the set

$I_{p}(A,B)=\{u_{0}\in H^{\frac{1}{2}}(\mathbb{T}^{3}):div u_{0}=0~and~(H1),~H(2),~(H3)~ are ~satisfied\}$, where

(H1)~~$\|Mu_{0}\|_{L^{2}(T^{2})}+\|MP(u_{F}\cdot\nabla u_{F})\|_{L^{1}(R^{+},L^{2}(T^{2}))}\leqslant A$;

(H2)~~$\|(Id-M)u_{0}\|_{B_{\infty,2}^{-1}}\leqslant A$;

(H3)~~$\|(Id-M)P(u_{F}\cdot\nabla u_{F})+Q(u_{2D},u_{F})\|_{L^{1}(R^{+},B_{p,2}^{-1+\frac{3}{p}})}\leqslant B,$\\
where we have noted $u_{F}=S(t)(Id-M)u^{0}$, $Mf(x_{1},x_{2})=\frac{1}{2\pi}\int_{0}^{2\pi}f(x_{1},x_{2},x_{3})dx_{3}$ and where $u_{2D}$ is a three component vector field defined
on $T^{2}$, satisfying the following two dimensional Navier-Stokes equation, in the case when the initial data is $\upsilon(0,x)=Mu^{0}(x)$, and the force is $f=-M\mathbb{P}(u_{F}\cdot \nabla)u_{F}$:
\begin{align*}(NS2D)\left\{
  \begin{array}{lll}
  \partial_{t}\upsilon-\Delta_{h} \upsilon+\mathbb{P}(\upsilon^{h}\cdot \nabla^{h})\upsilon=f,\\
  \upsilon(0,x)=\upsilon(0,x),\\
  \end{array}
  \right.
\end{align*}
where $\Delta_{h}$ denotes the horizontal Laplacian $\Delta_{h}=\partial_{1}^{2}+\partial_{2}^{2}$
 and where $\nabla^{h} = (\partial_{1},\partial_{2})$.\\\\
\begin{bf}Theorem D.\end{bf} \cite{Chemin1} Let $p\in(6,+\infty)$ be given. There is a constant $C_{0}>0$ such that the following
holds. Consider two positive real numbers A and B satisfying $$Bexp(C_{0}A^{2}(1+Alog(e+A))^{2})\leqslant C_{0}^{-1}.$$ Then for any vector field $u_{0}\in I_{p}(A,B)$, there is a unique, global solution $u$ to (NS) associated
with $u_{0}$, satisfying $$u\in C_{b}(R^{+};H^{\frac{1}{2}}(\mathbb{T}^{3}))\cap L^{2}(R^{+};H^{\frac{3}{2}}(\mathbb{T}^{3})).$$
\begin{bf}Theorem E.\end{bf} \cite{Liu} Let $\delta\in(0,1)$, $u_{0}=(u_{0}^{h},u_{0}^{3})\in H^{\frac{1}{2}}(\mathbb{R}^{3})\cap B_{2,1}^{0,\frac{1}{2}}(\mathbb{R}^{3})$ with $u_{0}^{h}$ belong to $L^{2}(\mathbb{R}^{3})\cap L^{\infty}(\mathbb{R}_{v};H^{-\delta}(\mathbb{R}^{2}))\cap L^{\infty}(\mathbb{R}_{v};H^{3}(\mathbb{R}^{2}))$. If we assume in addition
that $\partial_{3}u_{0}\in H^{-\frac{1}{2},0}$, then there exists a small enough positive constant $\varepsilon_{0}$ such that if $$\|\partial_{3}u_{0}\|_{H^{-\frac{1}{2},0}}^{2}exp\{C(A_{\delta}(u_{0}^{h})+B_{\delta}(u_{0}))\}\leqslant \varepsilon_{0},$$(NS) has a unique global solution $u\in C(\mathbb{R}^{+};H^{\frac{1}{2}})\cap L^{2}(\mathbb{R}^{+};H^{\frac{3}{2}})$, where
\begin{align*}A_{\delta}(u_{0}^{h})=&(\frac{\|\nabla_{h}u_{0}^{h}\|_{L_{v}^{\infty}(L_{h}^{2})}^{2}\|u_{0}^{h}\|_{L_{v}^{\infty}(B_{2,\infty}^{-\delta})_{h}}^{\frac{2}{\delta}}}{\|u_{0}^{h}\|_{L_{v}^{\infty}(L_{h}^{2})}^{\frac{2}{\delta}}}+\|u_{0}^{h}\|_{L_{v}^{\infty}(L_{h}^{2})}^{2})
exp\{C_{\delta}(1+\|u_{0}^{h}\|_{L_{v}^{\infty}(L_{h}^{2})}^{4})\},\end{align*}
\begin{align*}E_{\delta}(u_{0}^{h})=&\frac{\|u_{0}^{h}\|_{L_{v}^{\infty}(L_{h}^{2})}^{3}\|\nabla_{h}^{3}u_{0}^{h}\|_{L_{v}^{\infty}(L_{h}^{2})}^{\frac{1}{2}}}{\|\nabla_{h}u_{0}^{h}\|_{L_{v}^{\infty}(L_{h}^{2})}^{\frac{3}{2}}}+A_{\delta}(u_{0}^{h}),~~~~~~~~~\end{align*}
and\begin{align*}B_{\delta}(u_{0})=&\|u_{0}^{h}\|_{B_{2,1}^{0,\frac{1}{2}}}exp(CE_{\delta}(u_{0}^{h}))+\|u_{0}\|_{B_{2,1}^{0,\frac{1}{2}}}exp(\|u_{0}^{h}\|_{B_{2,1}^{0,\frac{1}{2}}}exp(CE_{\delta}(u_{0}^{h})))\end{align*}
are scaling invariant under the scaling transformation (2).\\\\
\begin{bf}Data availability statement.\end{bf}\\
The authors declare that the data supporting the findings of this study are available within
the paper.\\\\
\begin{bf}Conflicts of interest.\end{bf}\\
On behalf of all authors, the corresponding author states that there is no conflict of interest.\\\\
\begin{bf}Acknowledgement.\end{bf}\\
The research was supported by NSF of China (No. 11626218) and National key R$\&$D program of China (No. 2022YFA1005700).

 \label{}

%% The Appendices part is started with the command \appendix;
%% appendix sections are then done as normal sections
%% \appendix

%% \section{}
%% \label{}

%% References
%%
%% Following citation commands can be used in the body text:
%% Usage of \cite is as follows:
%%   \cite{key}         ==>>  [#]
%%   \cite[chap. 2]{key} ==>> [#, chap. 2]
%%

%% References with bibTeX database:

\bibliographystyle{elsarticle-num}
\bibliography{<your-bib-database>}

%% Authors are advised to submit their bibtex database files. They are
%% requested to list a bibtex style file in the manuscript if they do
%% not want to use elsarticle-num.bst.

%% References without bibTeX database:

% \begin{thebibliography}{00}

%% \bibitem must have the following form:
%%   \bibitem{key}...
%%

% \bibitem{}

% \end{thebibliography}

\end{document}